\documentclass[a4paper,leqno,draft]{article}
\usepackage{amssymb,amsmath,amsfonts,amsthm,
mathrsfs}

\newtheorem{theorem}{Theorem}[section]
\newtheorem{corollary}[theorem]{Corollary}
\newtheorem{lemma}[theorem]{Lemma}
\newtheorem{proposition}[theorem]{Proposition}
\newtheorem{definition}[theorem]{Definition}
\theoremstyle{definition}
\newtheorem{remark}[theorem]{Remark}
\newtheorem{example}[theorem]{Example}

\newcommand{\wt}[1]{\widetilde{#1}}

\newcommand{\Cinf}{\ensuremath{\mathcal{C}^\infty}}
\newcommand{\Cinfc}{\ensuremath{\mathcal{C}^\infty_{\text{c}}}}
\newcommand{\D}{\ensuremath{{\cal D}}}
\renewcommand{\S}{\mathscr{S}}
\newcommand{\E}{\ensuremath{{\cal E}}}

\renewcommand{\L}{\mathcal{L}}

\newcommand{\mb}[1]{\ensuremath{\mathbb{#1}}}
\newcommand{\N}{\mb{N}}

\newcommand{\R}{\mb{R}}
\newcommand{\C}{\mb{C}}

\newcommand{\G}{\ensuremath{{\cal G}}}
\newcommand{\Gt}{\ensuremath{{\cal G}_\tau}}
\newcommand{\GtS}{{\mathcal{G}}_{\tau,\S}}
\newcommand{\Gc}{\ensuremath{{\cal G}_\mathrm{c}}}
\newcommand{\Gcinf}{\ensuremath{{\cal G}^\infty_\mathrm{c}}}
\newcommand{\GS}{\G_{{\, }\atop{\hskip-4pt\scriptstyle\S}}\!}
\newcommand{\EM}{\ensuremath{{\cal E}_{M}}}

\newcommand{\Et}{\ensuremath{{\cal E}_{\tau}}}

\newcommand{\EMinf}{\ensuremath{{\cal E}^\infty_{M}}}

\newcommand{\Nt}{\ensuremath{{\cal N}_{\tau}}}
\newcommand{\Neg}{\mathcal{N}}
\newcommand{\NS}{\mathcal{N}_{\S}}

\newcommand{\Ginf}{\ensuremath{\G^\infty}}

\newcommand{\GSinf}{\G^\infty_{{\, }\atop{\hskip-3pt\scriptstyle\S}}}

\newcommand{\supp}{\mathrm{supp}}

\newfont{\bigmath}{cmr12 at 13pt}

\newfont{\grecomath}{cmmi12 at 15pt}

\newfont{\bl}{msbm10 scaled \magstep2}

\newcommand{\beq}{\begin{equation}}
\newcommand{\eeq}{\end{equation}}

\newcommand{\F}{\ensuremath{{\cal F}}}

\newcommand{\notmid}{\mid\kern-0.5em\not\kern0.5em}

\newcommand{\eps}{\varepsilon}

\newcommand{\Om}{\Omega}

\newcommand{\compl}[1]{{#1}^{\mathrm{c}}}

\newcommand{\val}{\mathrm{v}} 
\newcommand{\esp}{\mathrm{e}}
\newcommand{\M}{\mathcal{M}}

\newcommand{\mP}{\mathcal{P}}
\newcommand{\mQ}{\mathcal{Q}}
\newcommand{\mU}{\mathcal{U}}
\newcommand{\mF}{\mathcal{F}}
\newcommand{\mH}{\mathcal{H}}
\newcommand{\mO}{\mathcal{O}}
\newcommand{\bil}{\mathsf{b}}

\begin{document}
\title{\bf Topological structures in Colombeau algebras:
topological $\wt{\C}$-modules and duality theory}
\author{Claudia Garetto \footnote{Current address: Institut f\"ur Technische Mathematik, Geometrie und Bauniformatik,
Universit\"at Innsbruck, e-mail:\,\texttt{claudia@mat1.uibk.ac.at}}\\
Dipartimento di Matematica\\ Universit\`a di Torino,\\
via Carlo Alberto 10, 10123 Torino, Italia\\
\texttt{garettoc@dm.unito.it}\\
} 
\date{ }
\maketitle

\begin{abstract} We study modules over the ring $\wt{\C}$ of complex generalized numbers from a topological point of view, introducing the notions of $\wt{\C}$-linear topology and locally convex $\wt{\C}$-linear topology. In this context particular attention is given to completeness, continuity of $\wt{\C}$-linear maps and elements of duality theory for topological $\wt{\C}$-modules. As main examples we consider various Colombeau algebras of generalized functions.
\end{abstract}
 


\setcounter{section}{-1}
\section{Introduction}
Colombeau algebras of generalized functions have proved to be an analytically powerful tool in dealing with linear and nonlinear PDEs with highly singular coefficients \cite{Biagioni:90, BO:92b, BO:92, Colombeau:92, CHO:96, CO:90, GHKO:99, GH:03, HdH:01c, HdH:01, KO:00, LO:91, NPS:98, O:87, O:89, O:92, O:01, OW:99}. In the recent research on the
subject a variety of algebras of generalized functions \cite{BO:92, FGKS:01, Garetto:04, Grosser:01, GKOS:01, GH:04}
have been introduced in addition to the original construction by Colombeau \cite{Colombeau:84, Colombeau:85} and investigated in its algebraic and structural aspects as well as in analytic and applicative aspects. These investigations have produced a theory of point values in the Colombeau algebra $\G(\Om)$ and results of invertibility and positivity in the ring of constant generalized functions $\wt{\C}$ \cite{GKOS:01, OK:99, OPS:03} but also microlocal analysis in Colombeau algebras and regularity theory for generalized solutions to partial and (pseudo-) differential equations \cite{DPS:98, Garetto:04, GGO:03, GarH:04, Hoermann:99, GH:03, GH:04, HdH:01c, HdH:01, HK:01, HO:03, HOP:03}. Apart from some early and inspiring work by Biagioni, Pilipovi\'{c}, Scarpal\'{e}zos \cite{Biagioni:90, NPS:98, Scarpalezos:92, Scarpalezos:98, Scarpalezos:00}, topological questions have played a marginal role in the existing Colombeau literature. However, the recent papers on pseudodifferential operators acting on algebras of generalized functions \cite{Garetto:04, GGO:03, GarH:04} and a preliminary kernel theory introduced in \cite{GGO:03} motivate a renewed interest in topological issues, such as $\wt{\C}$-linear topologies on Colombeau algebras, $\wt{\C}$-linear continuous maps and duality theory.

This is the first of two papers devoted to a detailed topological investigation into algebras of generalized functions and succeeds to add a collection of original results to what is already known in the field. It develops a theory of topological $\wt{\C}$-modules and locally convex topological $\wt{\C}$-modules, which requires the introduction of $\wt{\C}$-versions of various concepts relating to topological and locally convex vector spaces. As a topic of particular interest, the foundations of duality theory are provided within this framework, dealing with the $\wt{\C}$-module $\L(\G,\wt{\C})$ of all $\wt{\C}$-linear and continuous functionals on $\G$. The second paper on topological structures in  Colombeau algebras \cite{Garetto:04c} will be focused on applications. Due to the fact that many algebras of generalized functions can be easily viewed as locally convex topological $\wt{\C}$-modules, we will be able to apply all the previous theoretical concepts and results to the \emph{topological dual of a Colombeau algebra}. This procedure together with the discussion of some relevant examples and continuous embeddings is a novelty in Colombeau theory.

We now describe the contents of the sections in more detail.

Section 1 serves to collect the basic topological notions which we will refer to in the course of the paper. Starting from the new notions of $\wt{\C}$- absorbent, balanced and convex subsets of a $\wt{\C}$-module $\G$, $\wt{\C}$-linear and locally convex $\wt{\C}$-linear topologies are introduced and described via their neighborhoods in Subsections 1.1 and 1.2 respectively. A characterization of locally convex topological $\wt{\C}$-modules is given, inspired by the analogous statements involving seminorms and locally convex vector spaces, making use of the concept of ultra-pseudo-seminorm. This turns out to be a useful technical tool in providing the classically expected results on separatedness and boundedness. In particular, the continuity of a $\wt{\C}$-linear map is expressed in terms of a uniform estimate between ultra-pseudo-seminorms. Inductive limits and strict inductive limits of locally convex topological $\wt{\C}$-modules are studied in Subsection 1.3. Finally Subsection 1.4 is concerned with completeness in topological $\wt{\C}$-modules. We pay particular attention to the relationships between completeness, strict inductive limit topology and initial topology in case of locally convex topological $\wt{\C}$-modules.

The theoretical core of the paper is Section 2 where we set the stage for the duality theory of topological $\wt{\C}$-modules. Using concepts as pairings of $\wt{\C}$-modules and polar sets we equip the dual $\L(\G,\wt{\C})$ with at least three locally convex $\wt{\C}$-linear topologies: the weak topology $\sigma(\L(\G,\wt{\C}),\G)$, the strong topology $\beta(\L(\G,\wt{\C}),\G)$ and the topology $\beta_b(\L(\G,\wt{\C}),\G)$ of uniform convergence on bounded subsets of $\G$. A theorem of completeness of the dual $\L(\G,\wt{\C})$ with respect to the strong topology as well as a $\wt{\C}$-linear formulation of the Banach-Steinhaus theorem are obtained under suitable hypotheses on $\G$. 

Section 3 investigates the properties of some interesting examples of locally convex topological $\wt{\C}$-modules and their topological duals. Inspired by \cite{Scarpalezos:98} Subsection 3.1 deals with the $\wt{\C}$-modules $\G_E$ of generalized functions based on the locally convex topological vector space $E$, showing how a separated locally convex $\wt{\C}$-linear topology may be defined, in terms of ultra-pseudo-seminorms, on $\G_E$ by means of the seminorms which topologize $E$. Well-known Colombeau algebras as $\wt{\C}$, $\G(\Om)$ \cite{GKOS:01} and $\GS(\R^n)$ \cite{Garetto:04, GGO:03} are recognized to be special cases. More sophisticated topological tools, as strict inductive limit topologies and initial topologies, are needed for the Colombeau algebra of compactly supported generalized functions $\Gc(\Om)$ \cite{GKOS:01} and the Colombeau algebra of tempered generalized functions $\Gt(\R^n)$ \cite{GKOS:01}. Finally, ultra-pseudo seminorms and norms fitted to measure the regularity of generalized functions are introduced, providing a topology for the Colombeau algebras $\Ginf(\Om)$, $\Ginf_{\rm{c}}(\Om)$, $\GSinf(\R^n)$ \cite{Garetto:04, GGO:03}. The continuity of $\wt{\C}$-linear maps of the form $T:\G_E\to\G$ is the topic of Subsection 3.2, while Subsection 3.3 is devoted to the topological dual $\L(\G_E,\wt{\C})$ when $E$ is a normed space.  In this particular case, an ultra-pseudo-norm modelled on the classical dual norm $\Vert\cdot\Vert_{E'}$ is defined on $\L(\G_E,\wt{\C})$ and a generalization of the Hahn-Banach theorem is given.
This result combined with a further adaptation of the Banach-Steinhaus theorem to the context of ultra-pseudo-normed $\wt{\C}$-modules allows to compare different $\wt{\C}$-linear topologies on $\L(\G_E,\wt{\C})$.

\section{Topological $\wt{\C}$-modules}
This section provides the required foundations of topology for $\wt{\C}$-modules. We begin with a collection of basic notions and definitions.

Let $\wt{\C}$ be the \emph{ring of complex generalized numbers} obtained factorizing $$\EM:=\{(u_\eps)_\eps\in\C^{(0,1]}:\, \exists N\in\N\quad |u_\eps|=O(\eps^{-N})\ \text{as}\ \eps\to 0\}$$ with respect to the ideal $$\Neg:=\{(u_\eps)_\eps\in\C^{(0,1]}:\, \forall q\in\N\quad |u_\eps|=O(\eps^{q})\ \text{as}\ \eps\to 0\}$$ (c.f. \cite{Colombeau:85, GKOS:01}).
$\wt{\C}$ is trivially a module over itself and it can be endowed with a structure of a topological ring. In order to explain this assertion, inspired by nonstandard analysis \cite{OT:98, TodorovW:00} and the previous work in this field \cite{Biagioni:90, NPS:98, Scarpalezos:92, Scarpalezos:98, Scarpalezos:00}, we introduce the function
\beq
\label{val_M}
\val: \EM\to(-\infty,+\infty]:(u_\eps)_\eps \to \sup\{b\in\R:\quad |u_\eps|=O(\eps^b)\, \text{as}\ \eps\to 0\}
\eeq
on $\EM$. It satisfies the following conditions:
\begin{itemize}
\item[(i)] $\val((u_\eps)_\eps)=+\infty$ if and only if $(u_\eps)_\eps\in\Neg$,
\item[(ii)] $\val((u_\eps)_\eps(v_\eps)_\eps)\ge \val((u_\eps)_\eps) +\val((v_\eps)_\eps)$, 
\item[(iii)] $\val((u_\eps)_\eps + (v_\eps)_\eps)\ge \min\{ \val((u_\eps)_\eps), \val((v_\eps)_\eps)\}$,
\end{itemize}
where $(ii)$ and $(iii)$ become equality if at least one or both terms are of the form $(c\eps^b)_\eps$, $c\in\C$, $b\in\R$, respectively. Note that if $(u_\eps-u'_\eps)_\eps\in\Neg$, $(i)$ combined with $(iii)$ yields $\val((u_\eps)_\eps)=\val((u'_\eps)_\eps)$. This means that we can use \eqref{val_M} to define the \emph{valuation} 
\beq
\label{val_C}
\val_{\wt{\C}}(u):=\val((u_\eps)_\eps)
\eeq
of the complex generalized number $u=[(u_\eps)_\eps]$, and that all the previous properties hold for the elements of $\wt{\C}$. Let now
\beq
\label{norm_C}
\vert\cdot\vert_\esp := \wt{\C}\to [0,+\infty): u\to \vert u\vert_\esp:=\esp^{-\val_{\wt{\C}}(u)}.
\eeq
The properties of the valuation on $\wt{\C}$ makes the coarsest topology on $\wt{\C}$ such that the map $\vert \cdot\vert_\esp$ is continuous compatible with the ring structure.
It is common in the already existing literature \cite{NPS:98, Scarpalezos:92, Scarpalezos:98, Scarpalezos:00} to use  the adjective ``sharp'' for such a topology. In this paper $\wt{\C}$ will always be endowed with its ``sharp topology''. 
Our investigation of the topological aspects of a $\wt{\C}$-module is mainly modeled on the classical approach to topological vector spaces and locally convex spaces suggested by many books on functional analysis \cite{Horvath:66, Robertson:80}. In particular it requires the adaptation of the algebraic notions of absorbent, balanced and convex subsets of a vector space, to the new context of $\wt{\C}$-modules.
\begin{definition}
\label{def_subsets}
A subset $A$ of a $\wt{\C}$-module $\G$ is $\wt{\C}$-absorbent if for all $u\in\G$ there exists $a\in\R$ such that $u\in[(\eps^b)_\eps]A$ for all $b\le a$.\\
$A\subseteq\G$ is $\wt{\C}$-balanced if $\lambda A\subseteq A$ for all $\lambda\in\wt{\C}$ with $\vert\lambda\vert_\esp\le 1$.\\
$A\subseteq\G$ is $\wt{\C}$-convex if $A+A\subseteq A$ and $[(\eps^b)_\eps]A\subseteq A$ for all $b\ge 0$.
\end{definition}
Note that if $A$ contains $0$ then it is $\wt{\C}$-convex if and only if $[(\eps^{b_1})_\eps]A+[(\eps^{b_2})_\eps]A\subseteq A$ for all $b_1,b_2\ge 0$. A subset $A$ which is both $\wt{\C}$-balanced and $\wt{\C}$-convex is called absolutely $\wt{\C}$-convex. In the case when $A$ is $\wt{\C}$-balanced the convexity is equivalent to the following statement: for all $\lambda,\mu\in\wt{\C}$ with $\max\{\vert\lambda\vert_\esp,\vert\mu\vert_\esp\}\le 1$, $\lambda A + \mu A\subseteq A$. The $\wt{\C}$-convexity cannot be considered as a generalization of the corresponding concept in vector spaces. In fact the only subset $A$ of $\C$ which is $\wt{\C}$-convex is the trivial set $\{0\}$.
\begin{remark}
\label{remark_grosser}
The definition of a $\wt{\C}$-balanced subset of $\G$ is inspired by the classical one concerning vector spaces and consists in replacing the absolute value in $\C$ with $|\cdot|_\esp$ in $\wt{\C}$. We may construct an analogy between a vector space $V$ and a $\wt{\C}$-module $\G$ by associating the sum in $V$ with the sum in $\G$ and the product $au$, $a>0$, $u\in V$, with $[(\eps^{-\log a})_\eps]u$ where $u\in\G$. In this way the concept of absorbent subset of $V$ is translated into the concept of $\wt{\C}$-absorbent subset of $\G$ and a convex cone (at 0) in $V$ corresponds to a $\wt{\C}$-convex subset of $\G$.
\end{remark}
In the sequel, we shall simply talk about absorbent, balanced or convex subset, omitting the prefix $\wt{\C}$, when we deal with $\wt{\C}$-modules. The reader should be aware that the words refer to Definition \ref{def_subsets} and not to the classical notions in this context. 
\subsection{Elementary properties of $\wt{\C}$-linear topologies}
We recall that a topology $\tau$ on a $\wt{\C}$-module $\G$ is said to be $\wt{\C}$-linear if the addition $\G\times\G\to\G:(u,v)\to u+v$ and the product $\wt{\C}\times\G\to \G:(\lambda,u)\to\lambda u$ are continuous.
A topological $\wt{\C}$-module $\G$ is a $\wt{\C}$-module with a $\wt{\C}$-linear topology. As an immediate consequence we have that for any $u_0\in\G$ and for any invertible $\lambda\in\wt{\C}$ the translation $\G\to\G:u\to u+u_0$ and the mapping $\G\to\G:u\to\lambda u$ are homeomorphisms of $\G$ into itself.
This means that if $\mathcal{U}$ is a base of neighborhoods of the origin that $\mathcal{U}+u_0$ is a base of neighborhoods of $u_0$ and if $U$ is a neighborhood of the origin so is $\lambda U$ for all invertible $\lambda\in\wt{\C}$. It is also clear that a $\wt{\C}$-linear map $T$ between topological $\wt{\C}$-modules is continuous if and only if it is continuous at the origin and then the set $\L(\G,\mH)$ of all continuous $\wt{\C}$-linear maps between the topological $\wt{\C}$-modules $\G$ and $\mH$ is a module on $\wt{\C}$.
\begin{proposition}
\label{prop_base}
Let $\G$ be a topological $\wt{\C}$-module and $\mathcal{U}$ be a base of neighborhoods of the origin. Then for each $U\in\mathcal{U}$,
\begin{itemize}
\item[(i)] $U$ is absorbent,
\item[(ii)] there exists $V\in\mathcal{U}$ with $V+V\subseteq U$,
\item[(iii)] there exists a balanced neighborhood of the origin $W$ such that $W\subseteq U$.
\end{itemize}
\end{proposition}
\begin{proof}
Fix $u\in\G$. The continuity of the product between elements of $\wt{\C}$ and elements of $\G$ guarantees for any $U\in\mU$ the existence of $\eta>0$ such that $\lambda u\in U$ for all $\lambda\in\wt{\C}$ with $\vert \lambda\vert_\esp\le\eta$. Hence $u\in[(\eps^b)_\eps]U$ provided $b\le\log\eta$. This shows that $U$ is absorbent.\\
The addition in $\G$ is continuous. Therefore, given $U\in\mU$ there exist $V_1,V_2\in\mU$ such that $V_1+V_2\subseteq U$ and a neighborhood $V\in\mU$ contained in $V_1\cap V_2$ which proves assertion $(ii)$.\\
Finally, since the product is continuous at $(0,0)$, there exist $\eta>0$ and $V\in\mU$ such that $\lambda V\subseteq U$ for $\vert\lambda\vert_\esp\le\eta$. Let $W=\cup_{\vert\lambda\vert_\esp\le\eta}\lambda V$. By construction $W$ is contained in $U$ and is a neighborhood of the origin since $[(\eps^{-\log\eta})_\eps]V\subseteq W$. Recalling that $\vert\lambda\mu\vert_\esp\le \vert\lambda\vert_\esp \vert\mu\vert_\esp$ for all complex generalized numbers $\lambda,\mu$, we conclude that $W$ is a balanced subset of $\G$.
\end{proof}
It follows from Proposition \ref{prop_base} that any topological $\wt{\C}$-module has a base of absorbent and balanced neighborhoods of the origin. As in the classical theory of topological vector spaces this fact ensures a useful characterization of separated topological $\wt{\C}$-modules. Further, if $\G$ is a topological $\wt{\C}$-module and $\mU$ a base of neighborhoods of the origin we have that $\G$ is separated if and only if $\displaystyle\cap_{U\in\mU}\, U=\{0\}$.
\begin{remark}
\label{remark_quotient_topology}
A particular example of a topological $\wt{\C}$-module is the quotient set $\G/M$, where $M$ is a $\wt{\C}$-submodule of the topological $\wt{\C}$-module $\G$, endowed with the quotient topology. In analogy with the theory of topological vector spaces, by the definition of a quotient $\wt{\C}$-module $\G/M$ and the previous considerations on separated modules we obtain that $\G/M$ is separated if and only if $M$ is closed in $\G$.
\end{remark}
Finally, having introduced a notion of absorbent set we can state the natural definition of a bounded subset of a topological $\wt{\C}$-module.
\begin{definition}
\label{def_bounded}
We say that a subset $A$ of a topological $\wt{\C}$-module $\G$ is bounded if it is absorbed by every neighborhood of the origin i.e. for all neighborhoods $U$ of the origin in $\G$ there exists $a\in\R$ such that $A\subseteq[(\eps^b)_\eps]U$ for all $b\le a$.
\end{definition}
A simple application of the definitions shows that any continuous $\wt{\C}$-linear map $T:\G\to \mH$ between topological $\wt{\C}$-modules is \emph{bounded}, in the sense that it maps bounded subsets of $\G$ into bounded subsets of $\mH$.

\subsection{Locally convex topological $\wt{\C}$-modules: ultra-pseudo-seminorms and continuity}
\begin{definition}
\label{def_convex_module}
A locally convex topological $\wt{\C}$-module is a topological $\wt{\C}$-module which has a base of $\wt{\C}$-convex neighborhoods of the origin.
\end{definition}
Proposition \ref{prop_base} shows that there exist bases of convex neighborhoods of the origin with additional properties.
\begin{proposition}
\label{prop_base_convex}
Every locally convex topological $\wt{\C}$-module $\G$ has a base of absolutely convex and absorbent neighborhoods of the origin.
\end{proposition}
\begin{proof}
Let $\mU$ be a base of convex neighborhoods of $0$ in $\G$. By Proposition \ref{prop_base}, for all $U\in\mU$ there exists a  a balanced neighborhood of the origin $W$ contained in $U$. Take the convex hull $W'$ of $W$ i.e. the set of all finite $\wt{\C}$-linear combinations of the form $[(\eps^{b_1})_\eps]w_1+[(\eps^{b_2})_\eps]w_2+...+[(\eps^{b_n})_\eps]w_n$ where $b_i\ge 0$ and $w_i\in W$. By construction $W'$ is an absolutely convex and absorbent neighborhood of $0$ and since $U$ is itself convex we have that $W'\subseteq U$.
\end{proof} 
We now want to deduce some more information on the topology of $\G$ from the nature of the neighborhoods. We begin with some preliminary definitions and results.
\begin{definition}
\label{def_ultra_pseudo}
Let $\G$ be a $\wt{\C}$-module. A \emph{valuation} on $\G$ is a function $\val:\G\to(-\infty,+\infty]$ such that
\begin{trivlist}
\item[(i)] $\val(0)=+\infty$,
\item[(ii)] $\val(\lambda u)\ge \val_{\wt{\C}}(\lambda)+\val(u)$ for all $\lambda\in\wt{\C}$, $u\in\G$,
\item[(ii)'] $\val(\lambda u)= \val_{\wt{\C}}(\lambda)+\val(u)$ for all $\lambda=[(c\eps^a)_\eps]$, $c\in\C$, $a\in\R$, $u\in\G$,
\item[(iii)] $\val(u+v)\ge\min\{\val(u),\val(v)\}$.
\end{trivlist}
An \emph{ultra-pseudo-seminorm} on $\G$ is a function $\mP:\G\to[0,+\infty)$ such that
\begin{trivlist}
\item[(i)] $\mP(0)=0$,
 \item[(ii)] $\mP(\lambda u)\le \vert\lambda\vert_\esp \mP(u)$ for all $\lambda\in\wt{\C}$, $u\in\G$,
\item[(ii)'] $\mP(\lambda u)= \vert\lambda\vert_\esp \mP(u)$ for all $\lambda=[(c\eps^a)_\eps]$, $c\in\C$, $a\in\R$, $u\in\G$,
\item[(iii)] $\mP(u+v)\le\max\{\mP(u),\mP(v)\}$.
\end{trivlist}
\end{definition}
The term valuation has here a slightly different meaning compared to the well-known concept introduced in nonstandard analysis and is deeply connected with the properties of $\G$ as a $\wt{\C}$-module. The reader should refer to \cite{LR:75, OT:98, Robinson:66, Robinson:73, TodorovW:00} for the original nonstardard approach and some related applications.

$\mP(u)=\esp^{-\val(u)}$ is a typical example of an ultra-pseudo-seminorm obtained by means of a valuation on $\G$. An \emph{ultra-pseudo-norm} is an ultra-pseudo-seminorm $\mP$ such that $\mP(u)=0$ implies $u=0$. $\vert\cdot\vert_\esp$ introduced in \eqref{norm_C} is an ultra-pseudo-norm on $\wt{\C}$. We now present an interesting example of a valuation on a $\wt{\C}$-module $\G$.
\begin{proposition}
\label{prop_gauge}
Let $A$ be an absolutely convex and absorbent subset of a $\wt{\C}$-module $\G$. Then 
\beq
\label{gauge}
\val_A(u):=\sup\{b\in\R:\, u\in[(\eps^b)_\eps]A\}
\eeq
is a valuation on $\G$. Moreover, for $\mP_A(u):=\esp^{-\val_A(u)}$ and $\eta>0$ the chain of inclusions
\beq
\label{chain}
\{u\in\G:\, \mP_A(u)<\eta\}\, \subseteq\, [(\eps^{-\log(\eta)})_\eps]A\, \subseteq\, \{ u\in\G:\, \mP_A(u)\le\eta\}
\eeq
holds.
\end{proposition}
We usually call $\mP_A$ the \emph{gauge} of $A$.
\begin{proof}
For each $u\in\G$ the set of real numbers $b$ such that $u\in [(\eps^b)]A$ is not empty. Hence $\val_A(u)$ is clearly a function from $\G$ into $(-\infty,+\infty]$. Since $A$ is balanced, $0$ belongs to $A$ and to every $[(\eps^b)_\eps]A$. Thus $\val_A(0)=+\infty$. Assume that $u\in[(\eps^b)]A$ for some $b\in\R$ and write 
\[
\lambda u = [(\eps^{b+\val_{\wt{\C}}(\lambda)})_\eps]\, \lambda\, [(\eps^{-\val_{\wt{\C}}(\lambda)})_\eps]\, [(\eps^{-b})_\eps]u,
\]
where $\lambda\in\wt{\C}\setminus 0$. From $\vert \lambda\, [(\eps^{-\val_{\wt{\C}}(\lambda)})_\eps]\vert_\esp = \vert\lambda\vert_\esp\, \esp^{\val_{\wt{\C}}(\lambda)}= 1$ and the fact that $A$ is a balanced subset of $\G$, we obtain that $\val_A(\lambda u)\ge \val_{\wt{\C}}(\lambda) + \val_A(u)$. In particular if $\lambda$ is of the form $[(c\eps^a)_\eps]$, $c\in\C\setminus 0$, $a\in\R$ and $\lambda u\in[(\eps^b)_\eps]A$, then 
\beq
\label{other_est}
u = [(\frac{1}{c}\eps^{-a})_\eps]\, [(\eps^b)_\eps]\, [(\eps^{-b})_\eps]\, [(c\eps^a)_\eps]u = [(\eps^{-a+b})_\eps]u'.
\eeq
Since $u'=[(\eps^{-b})_\eps][(\eps^a)_\eps]u\in A$, \eqref{other_est} leads to $\val_A(u)\ge -\val_{\wt{\C}}(\lambda) +\val_A(\lambda u)$ and shows $(ii)'$ in the definition a valuation.\\
Consider $u,v\in\G$. We know that there exist $b_1,b_2\in\R$ such that $u\in[(\eps^{b_1})_\eps]A$ and $v\in[(\eps^{b_2})_\eps]A$. The sum $u+v$ is an element of $[(\eps^{b_1}]A+[(\eps^{b_2})_\eps]A$ and we have that $u+v\in [(\eps^{b_1})_\eps](A+[(\eps^{b_2-b_1})_\eps]A)$. Let us assume $b_2-b_1\ge 0$ and recall that $A$ is convex. Hence $u+v\in[(\eps^{b_1})_\eps]A$ and $\val_A(u+v)\ge\min\{\val_A(u),\val_A(v)\}$.

Finally, in order to prove \eqref{chain} it is sufficient to observe that $\mP_A(u)<\eta$ implies $\val_A(u)>-\log(\eta)$ and $u\in[(\eps^{-\log(\eta)})_\eps]A$ while $u\in[(\eps^{-\log(\eta)})_\eps]A$ implies $\val_A(u)\ge -\log(\eta)$.
\end{proof}
\begin{theorem}
\leavevmode
\label{theorem_convex}
\begin{trivlist}
\item[(i)] Let $\{\mP_i\}_{i\in I}$ be a family of ultra-pseudo-seminorms on a $\wt{\C}$-module $\G$. The topology induced by $\{\mP_i\}$ on $\G$, i.e. the coarsest topology such that each ultra-pseudo-seminorm is continuous, induces the structure of locally convex topological $\wt{\C}$-module on $\G$.
\item[(ii)] In a locally convex topological $\wt{\C}$-module $\G$ the original topology is induced by the family of ultra-pseudo-seminorms $\{\mP_U\}_{U\in\mU}$, where $\mU$ is a base of absolutely convex and absorbent neighborhoods of the origin.
\end{trivlist}
\end{theorem}
\begin{proof}
\begin{trivlist}
\item[{\ }]
\item[(i)] From the properties $(ii)$ and $(iii)$ which characterize an ultra-pseudo-se\-mi\-norm it is clear that the coarsest topology such that the ultra-pseudo-seminorms $\{\mP_i\}_{i\in I}$ are continuous is $\wt{\C}$-linear on $\G$. A base of neighborhoods of the origin is given by all the finite intersections of sets of the form $\{u\in\G:\ \mP_i(u)\le\eta_i\}$ for $\eta_i>0$. Each $\{u\in\G:\ \mP_i(u)\le\eta_i\}$ is convex. In fact if $\mP_i(u_1)\le\eta_i$ and $\mP_i(u_2)\le\eta_i$ for all $b_1,b_2\ge 0$ we have that
\begin{multline*}
\mP_i([(\eps^b_1)_\eps]u_1+[(\eps^{b_2})_\eps]u_2)\le\max\{\mP_i([(\eps^{b_1})_\eps]u_1),\mP_i([(\eps^{b_2})_\eps]u_2)\}\\
=\max\{\esp^{-b_1}\mP_i(u_1),\esp^{-b_2}\mP_i(u_2)\}\le\eta_i.
\end{multline*}
Since a finite intersection of convex sets is still convex, $\G$ is a locally convex topological $\wt{\C}$-module.
\item[(ii)] Combining Proposition \ref{prop_gauge} with the previous considerations the topology induced by the family of ultra-pseudo-seminorms $\{\mP_U\}_{U\in\mU}$ is a locally convex $\wt{\C}$-linear topology on $\G$. \eqref{chain} relates the neighborhoods of the origin in this topology with the corresponding neighborhoods in the original topology on $\G$ and shows that the two topologies coincide.
\end{trivlist}
\end{proof}
Theorem \ref{theorem_convex} and the considerations after Proposition \ref{prop_base} lead to the following characterization of separated locally convex topological $\wt{\C}$-modules.
\begin{proposition}
\label{prop_T2_convex}
Let $\G$ be a locally convex topological ${\wt{\C}}$-module and $\{\mP_i\}_{i\in I}$ a family of continuous ultra-pseudo-seminorms which induces the topology of $\G$. $\G$ is separated if and only if for all $u\neq 0$ there exists $i\in I$ with $\mP_i(u)>0$.
\end{proposition}
\begin{example}
\label{example_quotient_convex}
If $\G$ is a locally convex topological $\wt{\C}$-module and $M$ is a $\wt{\C}$-submodule of $\G$ then $\G/M$ equipped with the quotient topology is locally convex itself. Note that if $\mQ$ is an ultra-pseudo-seminorm on $\G$ then, denoting the canonical projection of $\G$ on $\G /M$ by $\pi$,  
\[
\dot{\mQ}([u]):=\inf_{v\in\pi^{-1}([u])}\mQ(v)
\]
is a well-defined ultra-pseudo-seminorm on $\G/M$. Indeed, $\dot{\mQ}([0])=0$ and observing that for all invertible $\lambda$ in $\wt{\C}$, $v\in\pi^{-1}([\lambda u])$ if and only if $\lambda^{-1}v\in\pi^{-1}([u])$ we obtain the estimate 
\[
\dot{\mQ}(\lambda[u])=\inf_{v\in\pi^{-1}([\lambda u])}\mQ(v)\le\vert\lambda\vert_\esp\inf_{v\in\pi^{-1}([u])}\mQ(v)=\vert\lambda\vert_\esp\dot{\mQ}([u])
\]
which becomes an equality when $\lambda$ is of the form $[(c\eps^b)_\eps]\in\wt{\C}$. Finally, consider the sum $[u_1]+[u_2]$. If $\dot{\mQ}([u_1])<\dot{\mQ}([u_2])$ then for all $v_2\in\pi^{-1}([u_2])$ there exists $v_1\in\pi^{-1}([u_1])$ such that $\mQ(v_1)<\mQ(v_2)$ and this fact yields 
\[
\dot{\mQ}([u_1]+[u_2])\le\inf_{\substack{v_1\in\pi^{-1}([u_1])\\ v_2\in\pi^{-1}([u_2])}}\max\{\mQ(v_1),\mQ(v_2)\}\le \inf_{v_2\in\pi^{-1}([u_2])}\mQ(v_2)=\dot{\mQ}([u_2]).
\]
If $\dot{\mQ}([u_1])=\dot{\mQ}([u_2])$ then for all $v_2\in\pi^{-1}([u_2])$ and for all $\delta>0$ there exists $v_1\in\pi^{-1}([u_1])$ such that $\mQ(v_1)\le\mQ(v_2)+\delta$. It follows that  
\[
\dot{\mQ}([u_1]+[u_2])\le\inf_{\substack{v_1\in\pi^{-1}([u_1])\\ v_2\in\pi^{-1}([u_2])}}\max\{\mQ(v_1),\mQ(v_2)\}\le \inf_{v_2\in\pi^{-1}([u_2])}\mQ(v_2) +\delta
\]
and therefore $\dot{\mQ}([u_1]+[u_2])\le \dot{\mQ}([u_2])$.

The quotient topology $\tau$ on $\G /M$ is determined by the ultra-pseudo-seminorms $\{\dot{\mQ}\}_{\mQ}$, where $\mQ$ is an ultra-pseudo-seminorm continuous on $\G$. 
\end{example}
The ultra-pseudo-seminorms provide a useful tool for checking if a subset of a locally convex topological $\wt{\C}$-module is bounded. In the sequel let $(\G,\{\mP_i\}_{i\in I})$ be a locally convex topological $\wt{\C}$-module whose topology is determined by the family of ultra-pseudo-seminorms $\{\mP_i\}_{i\in I}$.
\begin{proposition}
\label{prop_bounded_sem}
Let $(\G,\{\mP_i\}_{i\in I})$ be a locally convex topological $\wt{\C}$-module. $A\subseteq\G$ is bounded if and only if for all $i\in I$ there exists a constant $C_i>0$ such that $\mP_i(u)\le C_i$ for all $u\in A$.
\end{proposition}
\begin{proof}
If $A\subseteq\G$ is bounded then for some $a_i\in\R$ it is contained in the set $[(\eps^{a_i})_\eps]\{u\in\G:\ \mP_i(u)\le 1\}$. This means that $\mP_i([(\eps^{-a_i})_\eps]u)\le 1$ for all $u\in A$ and the property $(ii)'$ which characterizes an ultra-pseudo-seminorm yields $\vert[(\eps^{-a_i})_\eps]\vert_\esp\mP_i(u)\le 1$. Thus $\mP_i(u)\le \esp^{-a_i}$ for every $u$ in $A$. Conversely, take a typical neighborhood of the origin of the form $U=\cap_{i\in I_0}\{u\in\G:\ \mP_i(u)\le \eta_i\}$ where $I_0$ is a finite subset of $I$. Again by the definition of an ultra-pseudo-seminorm and by $\mP_i(u)\le C_i$ on $A$ we have that $[(\eps^{-b})_\eps]A\subseteq U$ for all $b\le\log(\min_{i\in I_0}\eta_i)-\log(\max_{i\in I_0}C_i)$.
\end{proof} 
As in the classical theory of locally convex topological vector spaces an inspection of the neighborhoods of the origin gives some informations about ``metrizability'' and ``normability''. 
\begin{theorem}
\label{theorem_metric}
Let $\G$ be a separated locally convex topological $\wt{\C}$-module with a countable base of neighborhoods of the origin. Then its topology is induced by a metric $d$ invariant under translation.
\end{theorem}
\begin{proof}
Let $(U_n)_{n\in\N}$ be a countable base of neighborhoods of the origin in $\G$. By Proposition \ref{prop_base_convex} we may assume that each $U_n$ is absorbent and absolutely convex. We define
\[
f(u)=\sum_{n=0}^{\infty}2^{-n}\min\{\mP_{U_n}(u),1\},
\]
where $\mP_{U_n}$ is the gauge of $U_n$. By Proposition \ref{prop_T2_convex} $f(u)=0$ implies $u=0$, and by construction we have that $f(u)=f(-u)$, $f(u+v)\le f(u)+f(v)$ for all $u,v\in\G$. At this point, as in the proof of Theorem 4 in \cite[Chapter I]{Robertson:80} we obtain that $d(u,v):=f(u-v)$ is a distance invariant under translation which induces the original topology on $\G$. 
\end{proof}
The topology determined on a $\wt{\C}$-module $\G$ by an ultra-pseudo-norm $\mP$ is a separated and locally convex $\wt{\C}$-linear topology such that every set $\{u\in\G:\ \mP(u)\le \eta\}$ is bounded. This property characterizes the \emph{ultra-pseudo-normed} $\wt{\C}$-modules.
\begin{theorem}
\label{theorem_normed}
If $\G$ is a separated locally convex topological $\wt{\C}$-module and it has a bounded neighborhood of the origin, then the topology on $\G$ is induced by an ultra-pseudo-norm.
\end{theorem}
\begin{proof}
Let $V$ be an absorbent and absolutely convex neighborhood of the origin contained in a bounded neighborhood of the origin. Then $V$ is bounded, that is, for all neighborhoods $U$ of the origin in $\G$ there exists $a\in\R$ such that $V\subseteq[(\eps^b)_\eps]U$ for $b\le a$. This means that $[(\eps^{-b})_\eps]V\subseteq U$ and that $\{[(\eps^d)_\eps]V\}_{d\in\R}$ is a base of the neighborhoods of the origin in $\G$. Therefore, the gauge $\mP_V$ determines the topology on $\G$ which is separated. By Proposition \ref{prop_T2_convex}, $\mP_V$ is an ultra-pseudo-norm.
\end{proof}
We conclude this subsection with some continuity issues.  
\begin{theorem}
\label{theorem_sem}
Let $(\G,\{\mP_i\}_{i\in I})$ be a locally convex topological $\wt{\C}$-module. An ultra-pseudo-seminorm $\mQ$ on $\G$ is continuous if and only if it is continuous at the origin if and only if there exists a finite subset $I_0\subseteq I$ and a constant $C>0$ such that for all $u\in\G$
\beq
\label{est_gen_sem}
\mQ(u)\le C \max_{i\in I_0}\mP_i(u).
\eeq
\end{theorem}
\begin{proof}
Assume that $\mQ$ is continuous at the origin and take $u_0\in\G$, $u_0\neq 0$. For all $\delta>0$ there exists a finite subset $I_0\subseteq I$ and $\eta>0$ such that $\mQ(u)\le\delta$ if $\max_{i\in I_0}\mP_i(u)\le\eta$. Hence for all $u\in\G$ such that $\max_{i\in I_0}\mP_i(u-u_0)\le\eta$ we have that $\mQ(u-u_0)\le\delta$ and by definition of an ultra-pseudo-seminorm $|\mQ(u)-\mQ(u_0)|\le \mQ(u-u_0)$. This shows that $\mQ$ is continuous at $u_0\in\G$.\\
It is clear that if $\mQ$ satisfies \eqref{est_gen_sem} then it is continuous at the origin and consequently continuous on $\G$. Conversely if $\mQ$ is continuous at the origin as before there exists a finite subset $I_0\subseteq I$ and $\eta>0$ such that $\max_{i\in I_0}\mP_i(u)\le\eta$ implies $\mQ(u)\le 1$. We begin by observing that $\mQ(u)=0$ when $\max_{i\in I_0}\mP_i(u)=0$. In fact if $\mP_i(u)=0$ for all $i\in I_0$ then $$0=\vert[(\eps^b)_\eps]\vert_\esp \max_{i\in I_0}\mP_i(u) = \max_{i\in I_0}\mP_i([(\eps^b)_\eps]u)$$ and $\mQ([(\eps^b)_\eps]u)=\vert[(\eps^b)_\eps]\vert_\esp\mQ(u)=\esp^{-b}\mQ(u)\le 1$ for all $b\in\R$. So when the ultra-pseudo-seminorm $\max_{i\in I_0}\mP_i(u)$ is not zero we can write
\beq
\label{est_cont_Q}
\displaystyle
\mQ(v[(\eps^a)_\eps]) = \esp^{-a}\mQ(v),
\eeq
where $a=\log({\eta}/{\max_{i\in I_0}\mP_i(u)})$, $v=u[(\eps^{-a})_\eps]$ and by construction
$\mP_i(v)=\esp^a \mP_i(u)\le\eta$ for all $i\in I_0$. Combined with the continuity of $\mQ$ at the origin, \eqref{est_cont_Q} leads to \eqref{est_gen_sem} and completes the proof.
\end{proof}
Note that the composition of a $\wt{\C}$-linear map $T:\G\to\mH$ between $\wt{\C}$-modules with an ultra-pseudo-seminorm on $\mH$ gives an ultra-pseudo-seminorm on $\G$. Therefore, the following result concerning the continuity of $\wt{\C}$-linear maps between locally convex topological $\wt{\C}$-modules is a simple corollary of Theorem \ref{theorem_sem}.  
\begin{corollary}
\label{corollary_linear}
Let $(\G,\{\mP_i\}_{i\in I})$ and $(\mH,\{\mQ_j\}_{j\in J})$ be locally convex topological $\wt{\C}$-modules. A $\wt{\C}$-linear map $T:\G\to\mH$ is continuous if and only if it is continuous at the origin if and only if for all $j\in J$ there exists a finite subset $I_0\subseteq I$ and a constant $C>0$ such that for all $u\in\G$
\beq
\label{est_gen_lin}
\mQ_j(Tu) \le C \max_{i\in I_0}\mP_i(u).
\eeq
\end{corollary}

\subsection{Inductive limits and strict inductive limits of locally convex topological $\wt{\C}$-modules}
In this subsection we consider a family of locally convex topological $\wt{\C}$-modules $(\G_\gamma)_{\gamma\in\Gamma}$ and the $\wt{\C}$-module of all the finite $\wt{\C}$-linear combinations of elements of $\cup_{\gamma\in\Gamma}\G_\gamma$, denoted by span$(\cup_{\gamma\in\Gamma}\G_\gamma)$. We ask if the locally convex $\wt{\C}$-linear topologies $\tau_\gamma$ on $\G_\gamma$ can be pieced together to a locally convex $\wt{\C}$-linear topology $\tau$ on span$(\cup_{\gamma\in\Gamma}\G_\gamma)$. More generally we can start from a given $\wt{\C}$-module $\G$ which is spanned by the images under some $\wt{\C}$-linear maps $\iota_\gamma$ of the original $\G_\gamma$'s. 
\begin{theorem}
\label{theorem_inductive_limit}
Let $\G$ be a $\wt{\C}$-module, $(\G_\gamma)_{\gamma\in\Gamma}$ be a family of locally convex topological $\wt{\C}$-modules and $\iota_\gamma:\G_\gamma\to\G$ be a $\wt{\C}$-linear map so that $\G={\rm{span}}(\cup_{\gamma\in\Gamma}\iota_\gamma(\G_\gamma))$. Let
\[
\mU := \{U\subseteq \G\ \text{absolutely\ convex}: \forall\gamma\in\Gamma,\ \iota_\gamma^{-1}(U)\ \text{is a neighborhood of $0$ in}\ \G_\gamma\}.
\]
The topology $\tau$ induced by the gauges $\{\mP_U\}_{U\in\mU}$ is the finest $\wt{\C}$-linear topology with a base of absolutely convex neighborhoods of the origin such that each $\iota_\gamma$ is continuous.
\end{theorem}
With this topology $\G$ is called an \emph{inductive limit} of the locally convex topological $\wt{\C}$-modules $\G_\gamma$.
\begin{proof}
First we note that every $U\in\mU$ is absorbent. In fact, $\iota_\gamma^{-1}(U)$ is an absorbent neighborhood of $0$ in $\G_\gamma$ by Proposition \ref{prop_base} and then $U$ absorbs every element of $\iota_\gamma(\G_\gamma)$. Now when we take $u_1\in\G_{\gamma_1}$, $u_2\in\G_{\gamma_2}$, $\iota_{\gamma_1}(u_1)+\iota_{\gamma_2}(u_2)$ is absorbed by $U$ since we may write  
\begin{multline*}
\iota_{\gamma_1}(u_1)+\iota_{\gamma_2}(u_2)\in[(\eps^{b_1})_\eps]U+[(\eps^{b_2})_\eps]U\\=[(\eps^{\min\{b_1,b_2\}})_\eps]([(\eps^{b_1-\min\{b_1,b_2\}})_\eps]U+[(\eps^{b_2-\min\{b_1,b_2\}})_\eps]U),
\end{multline*}
for some $a_1,a_2\in\R$ and for all $b_1\le a_1$, $b_2\le a_2$, where, as observed after Definition \ref{def_subsets}, $[(\eps^{b_1-\min\{b_1,b_2\}})_\eps]U+[(\eps^{b_2-\min\{b_1,b_2\}})_\eps]U$ is contained in $U$. This means that $U$ is an absorbent subset of $\G$. By Proposition \ref{prop_gauge} and Theorem \ref{theorem_convex} the topology $\tau$ on $\G$ induced by the ultra-pseudo-seminorms $\{\mP_U\}_{U\in\mU}$ is a locally convex $\wt{\C}$-linear topology i.e. a $\wt{\C}$-linear topology with a base of absolutely convex neighborhoods of the origin. By definition of $\tau$ it is clear that every $\iota_\gamma:\G_\gamma\to (\G,\tau)$ is continuous. Assume now that $\tau'$ is another locally convex $\wt{\C}$-linear topology on $\G$ which makes each $\iota_\gamma$ continuous. $\tau$ is finer than $\tau'$ because if $U'$ is an absolutely convex neighborhood of $0$ for $\tau'$ then $\iota_\gamma^{-1}(U)$ is a neighborhood of $0$ in $\G_\gamma$ for all $\gamma\in\Gamma$.
\end{proof}
Continuity of $\wt{\C}$-linear maps between locally convex topological $\wt{\C}$-modules $\G$ and $\mH$ can easily be described when $\G$ has an inductive limit topology.  
\begin{proposition}
\label{prop_cont_induc}
Let $\G$ be the inductive limit of the locally convex topological $\wt{\C}$-modules $(\G_\gamma)_{\gamma\in\Gamma}$ and $\mH$ be a locally convex topological $\wt{\C}$-module. A $\wt{\C}$-linear map $T:\G\to\mH$ is continuous if and only if for each $\gamma\in\Gamma$ the composition $T\circ\iota_\gamma:\G_\gamma\to\mH$ is continuous.
\end{proposition}
\begin{proof}
The non-trivial assertion to prove is that $T$ is continuous if every $T\circ\iota_\gamma$ is continuous. By continuity at $0$, for every neighborhood $V$ of the origin in $\mH$, $\iota_\gamma^{-1}(T^{-1}(V))$ is a neighborhood of $0$ in $\G_\gamma$. Since we may choose $V$ absolutely convex and the $\wt{\C}$-linearity of $T$ guarantees that $T^{-1}(V)$ itself is absolutely convex in $\G$, the proof is complete.
\end{proof}
\begin{definition}
\label{def_strict}
Let $\G$ be a $\wt{\C}$-module and $(\G_n)_{n\in\N}$ be a sequence of $\wt{\C}$-sub\-mo\-du\-les of $\G$ such that $\G_n\subseteq \G_{n+1}$ for all $n\in\N$ and $\G=\cup_{n\in\N}\G_n$. Assume that $\G_n$ is equipped with a locally convex $\wt{\C}$-linear topology $\tau_n$ such that the topology induced by $\tau_{n+1}$ on $\G_n$ is $\tau_n$.\\
Then $\G$ endowed the inductive limit topology $\tau$ is called the \emph{strict inductive limit} of the sequence $(\G_n)_{n\in\N}$ of locally convex topological $\wt{\C}$-modules.
\end{definition}
\begin{proposition}
\label{prop_strict}
Let $\G$ be the strict inductive limit of the sequence of locally convex topological $\wt{\C}$-modules $(\G_n,\tau_n)_{n\in\N}$. The topology $\tau$ on $\G$ induces the original topology $\tau_n$ on each $\G_n$.
\end{proposition}
The proof of this proposition requires a technical lemma.
\begin{lemma}
\label{lemma1}
Let $\G$ be a locally convex topological $\wt{\C}$-module. Let $M$ be a $\wt{\C}$-submodule of $\G$ and $V$ be a convex neighborhood of the origin in $M$. Then there exists a convex neighborhood $W$ of $0$ in $\G$ such that $W\cap M=V$.
\end{lemma}
\begin{proof}
By definition of the induced topology on $M$ there exists a convex neighborhood $U$ of $0$ in $\G$ such that $U\cap M\subseteq V$. Let $W=U+V$. $W$ is the convex hull of $U\cup V$ since it can be written as $\{\sum_{i=1}^n[(\eps^{b_i})_\eps]u_i,\, n\in\N,\ b_i\ge 0,\ u_i\in U\cup V\}$, recalling the considerations after Definition \ref{def_subsets}. From $U\subseteq W$ we have that $W$ is a convex neighborhood of $0$ in $\G$ such that $V\subseteq W\cap M$.
It remains to prove the opposite inclusion. First, $w\in W$ is of the form $w=u+v$ for some $u\in U$ and $v\in V$. Therefore, if $w\in M$ we have that $u=w-v\in U\cap M$. Since $U\cap M\subseteq V$ we conclude that $u$ is an element of $V$. This leads to $W\cap M\subseteq V$.
\end{proof}
\begin{remark}
\label{rem_lemma1}
Note that if $U$ and $V$ are both convex and balanced then $U+V$ is the absolutely convex hull of $U\cup V$ i.e. the set of all finite sums $\sum_{i=1}^n\lambda_i u_i$ where $u_i\in U\cup V$, $\lambda_i\in\wt{\C}$ and $\max_{i=1,...,n}\vert\lambda_i\vert_\esp\le 1$.
\end{remark}  
\begin{proof}[Proof of Proposition \ref{prop_strict}]
Denoting the topology induced by $\tau$ on $\G_n$ by $\tau'_n$, it is clear that $\tau'_n$ is coarser than $\tau_n$. It remains to prove that any absolutely convex neighborhood $V_n$ of the origin in the topology $\tau_n$ is obtained as the intersection of a neighborhood of $0$ in $\G$ with the $\wt{\C}$-module $\G_n$. Lemma \ref{lemma1} and Remark \ref{rem_lemma1} allow us to construct a sequence $(V_{n+p})_{p\in\N}$ such that $V_{n+p}$ is an absolutely convex neighborhood of the origin in $\G_{n+p}$ for $\tau_{n+p}$, $V_{n+p}\subseteq V_{n+p+1}$ and $V_{n+p}\cap\G_n=V_n$ for all $p$. In conclusion, $V=\cup_{p\in\N}V_{n+p}$ is a neighborhood of the origin in $\G$ such that $V\cap\G_n=V_n$. 
\end{proof}
The following statements concerning separated $\wt{\C}$-modules and the closedness of $\G_n$ in $\G$ are immediate consequences of Proposition \ref{prop_strict}. We refer to \cite[Chapter 2, Section 12, Cor. 1,2]{Horvath:66} for a proof.
\begin{corollary}
\label{corollary_strict}
Under the hypotheses of Proposition \ref{prop_strict}, $\G$ is separated if each $\G_n$ is separated.
\end{corollary}
\begin{corollary}
\label{corollary_strict_closed}
Under the hypotheses of Proposition \ref{prop_strict}, if each $\G_n$ is closed in $\G_{n+1}$ for the topology $\tau_{n+1}$ then $\G_n$ is closed in $\G$ for $\tau$.
\end{corollary}  
We conclude the collection of results involving strict inductive limits of locally convex topological $\wt{\C}$-modules by characterizing bounded subsets.
\begin{theorem}
\label{theorem_bounded_strict}
Let $(\G,\tau)$ be the strict inductive limit of the sequence of locally convex topological $\wt{\C}$-modules $(\G_n,\tau_n)_{n\in\N}$. Assume in addition that each $\G_n$ is closed in $\G_{n+1}$ with respect to $\tau_{n+1}$. Then $A\subseteq\G$ is bounded if and only if $A$ is contained in some $\G_n$ and bounded there.
\end{theorem}
The proof of Theorem \ref{theorem_bounded_strict} requires some preliminary lemmas.
\begin{lemma}
\label{lemma2}
A set $A$ in a topological $\wt{\C}$-module $\G$ is bounded if and only if for all sequences $(u_n)_n$ of elements of $A$ and all sequences $(\lambda_n)_n$ in $\wt{\C}$ converging to $0$, the sequence $(\lambda_n u_n)_n$ tends to $0$ in $\G$.
\end{lemma}
\begin{proof}
Let $A$ be a bounded subset of $\G$ and $V$ be a balanced neighborhood of the origin. Since $A\subseteq[(\eps^a)_\eps]V$ for some $a\in\R$ we have that $\mP_V(u)\le \esp^{-a}$ on $A$. As shown in the proof of Proposition \ref{prop_gauge}, the estimate $\mP_V(\lambda u)\le\vert\lambda\vert_\esp\mP_V(u)$ holds for all $\lambda\in\wt{\C}$ and $u\in\G$. Therefore $\mP_V(\lambda_n u_n)\le\vert\lambda_n\vert_\esp\mP_V(u_n)\le\vert\lambda_n\vert_\esp\esp^{-a}$ and $\lambda_n\to 0$ in $\wt{\C}$ yields $\mP_V(\lambda_n u_n)<1$ for $n$ larger than some $N\in\N$. As a consequence, $\lambda_n u_n\in V$ if $n\ge N$ and $\lambda_n u_n$ is convergent to $0$ in $\G$.\\
Suppose now that all the sequences $(\lambda_n u_n)_n$, where $(u_n)_n\subseteq A$ and $\lambda_n\to 0$ in $\wt{\C}$, tend to $0$ in $\G$. Then $A$ is necessarily bounded. In fact if $A$ is not bounded there exists a balanced neighborhood of the origin $U$ and a sequence $b_n\to -\infty$ such that $A\cap (\G\setminus[(\eps^{b_n})_\eps]U)\neq\emptyset$. Choosing $u_n\in A\cap (\G\setminus[(\eps^{b_n})_\eps]U)$, the sequence $[(\eps^{-b_n})_\eps]$ goes to $0$ in $\wt{\C}$ but $[(\eps^{-b_n})_\eps]u_n$ is not convergent to $0$ in $\G$ since $[(\eps^{-b_n})_\eps]u_n\not\in U$ for all $n\in\N$. This contradicts our hypothesis.
\end{proof}
\begin{lemma}
\label{lemma1'}
Under the assumptions of Lemma \ref{lemma1}, if $M$ is closed then for every $u_0\not\in M$ there exists a convex neighborhood $W_0$ of $0$ in $\G$ such that $W_0\cap M= V$ and $u_0\not\in W_0$. 
\end{lemma}
\begin{proof}
If $M$ is closed then by Remark \ref{remark_quotient_topology} and Example \ref{example_quotient_convex} $\G/M$ is a separated locally convex topological $\wt{\C}$-module. This implies that there exists a convex neighborhood $U_0$ of $0$ in $\G$ such that $[u_0]\not\in\pi(U_0)$. Hence $(u_0+M)\cap U_0=\emptyset$. By Lemma \ref{lemma1} there exists a convex neighborhood $W$ of $0$ in $\G$ such that $W\cap M=V$. Therefore taking $W\cap U_0$, we can state that there exists a convex neighborhood $U'_0$ of $0$ in $\G$ such that $(u_0+M)\cap U'_0=\emptyset$ and $U'_0\cap M\subseteq V$. The same reasoning as in Lemma \ref{lemma1} combined with $(u_0+M)\cap U'_0=\emptyset$ shows that $W_0=U'_0+V$ is a convex neighborhood of $0$ in $\G$ such that $W_0\cap M=V$ and $u_0\not\in W_0$.
\end{proof} 
Note that if we choose $V$ and $U_0$ absolutely convex then by Remark \ref{rem_lemma1} we obtain that $W_0$ is an absolutely convex neighborhood of the origin in $\G$.
\begin{proof}[Proof of Theorem \ref{theorem_bounded_strict}]
If $A\subseteq \G_n$ is bounded for the topology $\tau_n$ then the continuity of the embedding of $(\G_n,\tau_n)$ into $(\G,\tau)$ guarantees that $A$ is bounded in $\G$.\\
Suppose now that $A$ is not contained in any $\wt{\C}$-module $\G_n$ and choose a sequence of elements $u_n\in A\cap(\G\setminus\G_n)$. There exists a strictly increasing sequence $(n_k)_k$ of natural numbers and a subsequence $(v_k)_k$ of $(u_n)_n$ such that $v_k\in\G_{n_{k+1}}\setminus \G_{n_k}$. By Lemma \ref{lemma1'} we can construct an increasing sequence $(V_k)_k$ of absolutely convex sets such that $V_k$ is a neighborhood of $0$ in $\G_{n_k}$, $V_{k+1}\cap \G_{n_k}=V_k$ and $[(\eps^{k})_\eps]v_k\not\in V_{k+1}$. As in the proof of Proposition \ref{prop_strict}, $V=\cup_{k\in\N}V_k$ is a neighborhood of the origin in $\G$ which does not contain $[(\eps^k)_\eps]v_k$ for any $k\in\N$. Then $[(\eps^k)_\eps]\to 0$ in $\wt{\C}$ but the sequence $([(\eps^k)_\eps]v_k)_k$ is not convergent to $0$ in $\G$. By Lemma \ref{lemma2} it follows that $A$ cannot be bounded in $\G$.\\
Finally, by Proposition \ref{prop_strict} it is clear that if $A$ is contained in some $\G_n$ and bounded in $\G$ it has to be bounded in $\G_n$ as well. 
\end{proof}
Every sequence $(u_n)_n$ in $\G$ which is tending to $0$ is an example of bounded set in $\G$. In fact for each absolutely convex neighborhood $U$ of the origin, there exists $N\in\N$ such that $u_n\in U$ for all $n\ge N$, and noting that $[(\eps^{b_1})_\eps]U\subseteq [(\eps^{b_2})_\eps]U$ if $b_1\ge b_2$, there exists $a\in\R$ such that $u_n\in[(\eps^b)_\eps]U$ for all $n\in\N$ and $b\le a$. At this point recalling Proposition \ref{prop_strict} it is immediate to prove the following corollary of Theorem \ref{theorem_bounded_strict}.
\begin{corollary}
\label{cor_cont_strict}
Under the assumptions of Theorem \ref{theorem_bounded_strict} a sequence $(u_n)_n$ is convergent to $0$ in $\G$ if and only if it is contained in some $\G_n$ and convergent to $0$ there.
\end{corollary}

\subsection{Completeness}
In this subsection we adapt the theory of complete topological vector spaces \cite{Horvath:66} to the context of $\wt{\C}$-modules. We say that a subset $A$ of a topological $\wt{\C}$-module is \emph{complete} if every Cauchy filter on $A$ converges to some point of $A$ and that a topological $\wt{\C}$-module $\G$ is \emph{quasi-complete} if every bounded closed subset is complete.
\begin{remark}
\label{remark_complete}
A topological $\wt{\C}$-module $\G$ is a uniform space \cite{Bourbaki:66, SchaeferW:99} and hence the following properties hold which will be used repeatedly later.  

Let $A$ be a subset of $\G$. Any filter $\mF$ on $A$ convergent to some $u\in \G$ is a Cauchy filter and if $u\in\G$ adheres to a Cauchy filter $\mO$ on $A$ then $\mO$ converges to $u$. Moreover a complete subset of a separated topological $\wt{\C}$-module is closed and if $A\subseteq \G$ is complete every closed subset of $A$ is complete itself. Finally in a metrizable topological $\wt{\C}$-module $\G$ a subset $A$ is complete if and only if every Cauchy sequence of points of $A$ converges to some point of $A$.
\end{remark}
Note that even if $\G$ is only quasi-complete, every Cauchy sequence $(u_n)_n\subseteq\G$ is convergent. First of all since $(u_n)_n$ is a Cauchy sequence the set $U:=\{u_n ,\ n\in\N\}$ is bounded in $\G$. We recall that for all neighborhoods $V$ of $0$ in a topological $\wt{\C}$-module we can find a balanced neighborhood $W$ of $0$ such that $W+W\subseteq V$. This means that $\overline{W}\subseteq V$ and therefore the closure of a bounded subset of a topological $\wt{\C}$-module is still bounded. Then in our case $\overline{U}$ is closed and bounded in $\G$ and by the quasi-completeness of $\G$ the sequence $(u_n)_n\subseteq \overline{U}$ is convergent.
 
A locally convex topological $\wt{\C}$-module which is metrizable and complete is called a \emph{Fr\'echet $\wt{\C}$-module}. As a straightforward application of Remark \ref{remark_complete} we show that $\wt{\C}$ is complete. The proof of this result is essentially due to Scarpal\'{e}zos \cite[Proposition 2.1]{Scarpalezos:98}.
\begin{proposition}
\label{prop_C_complete}
$\wt{\C}$ with the topology given by the ultra-pseudo-norm $\vert\cdot\vert_\esp$ is complete.
\end{proposition}
\begin{proof}
By Remark \ref{remark_complete} it is sufficient to prove that every Cauchy sequence $(u_n)_n$ in $\wt{\C}$ is convergent. We know that for every $\eta>0$ there exists $N\in\N$ such that for all $m,p\ge N$, $\vert u_m-u_p\vert_\esp\le\eta$. Considering representatives and the valuation on $\wt{\C}$ defined via $\val:\E_M\to(-\infty,+\infty]$ in \eqref{val_M}, we can extract a subsequence $(u_{n_k})_k$ such that $\val((u_{n_{k+1},\eps}-u_{{n_k},\eps})_\eps)> k$ for all $k\in\N$. This means that we can find $\eps_k\searrow 0$, $\eps_k\le 1/2^k$ such that $|u_{n_{k+1},\eps}-u_{n_k,\eps}|\le \eps^k$ on $(0,\eps_k)$. Let  
\[
h_{k,\eps}=\begin{cases} u_{n_{k+1},\eps}-u_{n_k,\eps}\  & \eps\in(0,\eps_k),\\
0\ & \eps\in[\eps_k,1].
\end{cases}
\]
$(h_{k,\eps})_\eps\in\E_M$ since $|h_{k,\eps}|\le\eps^k$ on the interval $(0,1]$. Moreover the sum  $u_\eps:=u_{n_0,\eps}+\sum_{k=0}^\infty h_{k,\eps}$ is locally finite and by
\[
|u_{\eps}|\le |u_{n_0,\eps}|+ \sum_{k=0}^\infty|h_{k,\eps}|\le |u_{n_0,\eps}| + \sum_{k=0}^\infty\frac{1}{2^k}
\]
it defines the representative of a complex generalized number $u=[(u_\eps)_\eps]$. The sequence $u_{n_k}$ tends to $u$ in $\wt{\C}$. In fact, for all $\overline{k}\ge 1$ the estimate
\[
|u_{n_{\overline{k}},\eps}-u_\eps|=\big|-\sum_{k=\overline{k}}^{\infty}h_{k,\eps}\big|\le \sum_{k=\overline{k}}^{\infty}\eps^{k-1}\eps_k\le \eps^{\overline{k}-1}\sum_{k=\overline{k}}^{\infty}\frac{1}{2^k},
\]
valid on the interval $(0,\eps_{\overline{k}-1})$, yields $\val((u_{n_k,\eps}-u_\eps)_\eps)\to +\infty$. Thus $(u_n)_n$ is a Cauchy sequence with a convergent subsequence and it converges to the same point $u\in\wt{\C}$.
\end{proof}
It is possible to decide if a strict inductive limit of locally convex topological $\wt{\C}$-modules is complete by looking at the terms $\G_n$ of the sequence which defines it.
\begin{theorem}
\label{theorem_strict_complete}
Let $(\G,\tau)$ be the strict inductive limit of the sequence of locally convex topological $\wt{\C}$-modules $(\G_n,\tau_n)_{n\in\N}$ where $\G_n$ is assumed to be closed in $\G_{n+1}$ for $\tau_{n+1}$. Then $\G$ is complete if and only if all the $\G_n$ are complete.
\end{theorem}
Before proving this theorem we present a technical lemma which will turn out to be useful later on as well. 
\begin{lemma}
\label{lemma_filter}
Let $\mF$ be a Cauchy filter on the strict inductive limit $(\G,\tau)$ of Theorem \ref{theorem_strict_complete} and $\mO$ be the Cauchy filter whose base is formed by all the sets $M+V$ where $M$ runs through $\mF$ and $V$ through the filter of neighborhoods of the origin in $\G$. Then there exists an integer $n$ such that $\mO$ induces a Cauchy filter on $\G_n$.
\end{lemma}
\begin{proof}
If there exists $n\in\N$ such that $(M+V)\cap \G_n\neq\emptyset$ for all $M\in\mF$ and neighborhoods $V$ of the origin in $\G$ then by Proposition \ref{prop_strict} the lemma is proven. We assume therefore that this is not the case, i.e. that for all $n\in\N$ there exist $M_n\in\mF$ and a neighborhood $V_n$ of the origin in $\G$ such that $(M_n+V_n)\cap\G_n=\emptyset$. In addition we may assume that $M_n-M_n\subseteq V_n$ and that $(V_n)_n$ is a decreasing sequence of absolutely convex neighborhoods. Consider the absolutely convex hull $W$ of $\cup_{n\in\N}(V_n\cap\G_n)$. Since every $V_n$ is absolutely convex it coincides with the set of all finite sums of elements of $\cup_{n\in\N}(V_n\cap\G_n)$ and by construction it is a neighborhood of the origin in $\G$. We want to prove that no $Q\in\mF$ has the property $Q-Q\subseteq W$. For this purpose we take $W_n:=V_0\cap\G_0+V_1\cap\G_1+\dots +V_{n-1}\cap\G_{n-1}+V_n$ which is the absolutely convex hull of $(\cup_{i\le n-1}V_i\cap\G_i)\cup V_n$. $W_n$ is a neighborhood of the origin in $\G$ and $W\subseteq W_n$ for all $n$. Since $\mF$ is a Cauchy filter there exists $P_n\in\mF$ such that $P_n-P_n\subseteq W_n$. This implies $(P_n+W_n)\cap\G_n=\emptyset$. In fact for $u_0\in P_n\cap M_n$ we have that $P_n\subseteq u_0+W_n$ and as a consequence every element $y$ of $P_n$ has the form $y=u_0+\sum_{i=0}^n v_i$ where $v_i\in V_i\cap\G_i$ if $i=0,...,n-1$ and $v_n\in V_n$. At this point $z\in P_n+W_n$ may be written as $z=u_0+\sum_{i=0}^{n-1}(v_i+v'_i) + v_n+v'_n$ with $v_i, v'_i\in V_i\cap\G_n$ for $i=0,...,n-1$ and $v_n,v'_n\in V_n$. Note that since $V_n$ is convex $v_n+v'_n\in V_n$ and that $\sum_{i=0}^{n-1}(v_i+v'_i)\in \G_n$. By $(M_n+V_n)\cap\G_n=\emptyset$ it follows that $z\not\in\G_n$.\\
Finally suppose that there exists $Q\in\mF$ such that $Q-Q\subseteq W$ and that $y_0\in Q$. Then $y_0\in\G_n$ for some $n$ and $Q\cap P_n=\emptyset$ which contradicts the hypothesis that $\mF$ is a filter. Indeed by construction of $P_n$ and $W_n$ if $y\in P_n$ then $y_0-y\in\compl{W_n}\subseteq\compl{W}$. Hence $y$ does not belong to $Q$.  
\end{proof}
\begin{proof}[Proof of Theorem \ref{theorem_strict_complete}]
If $\G$ is complete, recalling that by Corollary \ref{corollary_strict_closed} every $\G_n$ is closed in $(\G,\tau)$, by Remark \ref{remark_complete} every $\G_n$ is complete.
Conversely assume that each $(\G_n,\tau_n)$ is complete and take a Cauchy filter $\mF$ on $\G$. The filter $\mO$ constructed in Lemma \ref{lemma_filter} induces a Cauchy filter $\mO_n$ on some $\G_n$. Hence $\mO_n$ converges to some $u\in\G_n$ and, since by Proposition \ref{prop_strict} $\tau_n$ is the topology induced by $\tau$ on $\G_n$, $u$ adheres to the filter $\mO$. Consequently $\mO$ converges to $u$ and the same conclusion holds for $\mF$, since it is finer than $\mO$.
\end{proof}
We finally consider a family of topological $\wt{\C}$-modules $(\G_\gamma)_{\gamma\in\Gamma}$ and a $\wt{\C}$-module $\G$ such that for each $\gamma\in \Gamma$ there exists a $\wt{\C}$-linear map $\iota_\gamma:\G\to \G_\gamma$. The \emph{initial topology} on $\G$ is the coarsest topology such that each $\iota_\gamma$ is continuous. By the $\wt{\C}$-linearity of $\iota_\gamma$ we have that such a topology is $\wt{\C}$-linear and a base of neighborhoods of the origin is given by all the finite intersections $\iota_{\gamma_1}^{-1}(U_1)\cap\iota_{\gamma_2}^{-1}(U_2)...\cap\iota_{\gamma_n}^{-1}(U_n)$ where $U_{i}$, $i=1,2,...,n$, is a neighborhood of $0$ in $\G_{\gamma_i}$. In particular if the $\G_\gamma$ are locally convex topological $\wt{\C}$-modules with ultra-pseudo-seminorms $\{\mP_{j,\gamma}\}_{j\in J_\gamma}$ then the initial topology on $\G$ is determined by the family of ultra-pseudo-seminorms $\{\mP_{j,\gamma}\circ\iota_\gamma\}_{j\in J_\gamma, \gamma\in\Gamma}$. Let now $I$ be an ordered set of indices and $(\G_i)_{i\in I}$ be a family of topological $\wt{\C}$-modules such that $\G_j\subseteq\G_i$ if $j\ge i$. The intersection $\G:=\cap_{i\in I}\G_i$ is naturally endowed with the initial topology defined by $(\G_i)_{i\in I}$ and the injections $\G\to\G_i$. Adapting the reasoning of Proposition 3 and the corresponding corollary in \cite[Chapter 2, Section 11]{Horvath:66} to the context of topological $\wt{\C}$-modules, we prove that the completeness of each $\G_i$ may be transferred to the intersection $\G$ under suitable hypotheses.
\begin{proposition}
\label{prop_complete_initial}
Let $(\G_i)_{i\in I}$ be a family of separated topological $\wt{\C}$-modules where the index set is ordered. Suppose that for $i\le j$, $\G_j$ is a $\wt{\C}$-submodule of $\G_i$ and the topology on $\G_j$ is finer than the topology induced by $\G_i$ on $\G_j$. Let $\G=\cap_{i\in I}\G_i$ be equipped with the initial topology for the injections $\G\to \G_i$. If the $\G_i$ are complete then $\G$ is complete.
\end{proposition}

\section{Duality theory for topological $\wt{\C}$-modules}
This section is devoted to the dual of a topological $\wt{\C}$-module $\G$ i.e. the $\wt{\C}$-module $\L(\G,\wt{\C})$ of all $\wt{\C}$-linear and continuous maps on $\G$ with values in $\wt{\C}$. We present different ways of endowing $\L(\G,\wt{\C})$ with a $\wt{\C}$-linear topology and deal with related topics as pairings of $\wt{\C}$-modules, weak topologies, polar sets and polar topologies. 
\begin{definition}
\label{def_pairing}
Let $\G$ and $\mH$ be two $\wt{\C}$-modules. If a $\wt{\C}$-bilinear form $\bil:\G\times\mH\to\wt{\C}:(u,v)\to\bil(u,v)$ is given we say that $\G$ and $\mH$ form a pairing with respect to $\bil$. The pairing separates points of $\G$ if for all $u\neq 0$ in $\G$ there exists $v\in\mH$ such that $\bil(u,v)\neq 0$. Analogously it separates points of $\mH$ if for all $v\neq 0$ in $\mH$ there exists $u\in\G$ such that $\bil(u,v)\neq 0$. The pairing is separated if it separates points of both $\G$ and $\mH$.
\end{definition}
A pairing $(\G,\mH,\bil)$ defines a topology on each involved $\wt{\C}$-module via the $\wt{\C}$-bilinear form $\bil$. The \emph{weak topology} on $\G$ is the coarsest topology $\sigma(\G,\mH)$ on $\G$ such that each map $\bil(\cdot,v):\G\to\wt{\C}:u\to\bil(u,v)$, for $v$ varying in $\mH$, is continuous. Every $\bil(\cdot,v)$ is $\wt{\C}$-linear and continuous if and only if the ultra-pseudo-seminorm $\mP_v:\G\to[0,\infty):u\to\vert\bil(u,v)\vert_\esp$ is continuous. Hence $\sigma(\G,\mH)$ is the topology induced by the family of ultra-pseudo-seminorms $\{\mP_v\}_{v\in\mH}$ and by Theorem \ref{theorem_convex} it provides the structure of a locally convex topological $\wt{\C}$-module on $\G$. Obviously the same holds for $(\mH, \sigma(\mH,\G))$.

Combining Definition \ref{def_pairing} with Proposition \ref{prop_T2_convex} we obtain that the pairing $(\G,\mH,\bil)$ separates points of $\G$ if and only if $\sigma(\G,\mH)$ is a Hausdorff topology. Any $\wt{\C}$-module $\G$ with its algebraic dual $L(\G,\wt{\C})$ and any topological $\wt{\C}$-module $\G$ with its topological dual $\L(\G,\wt{\C})$ forms a pairing via the canonical $\wt{\C}$-bilinear map $\langle u,T\rangle=T(u)$. By the previous considerations the topologies $\sigma(L(\G,\wt{\C}),\G)$ and $\sigma(\L(\G,\wt{\C}),\G)$ are separated.
\begin{proposition}
\label{complete_algebraic_dual}
Let $\G$ be a $\wt{\C}$-module. $L(\G,\wt{\C})$ is complete for the weak topology $\sigma(L(\G,\wt{\C}),\G)$.
\end{proposition}
\begin{proof}
Let $\mF$ be a Cauchy filter on $L(\G,\wt{\C})$. For all $u\in\G$, $\mF_u$, the filter having as a base the family $\{X_u\}_{X\in\mF}$ where $X_u:=\{T(u):\, T\in X\}$, is a Cauchy filter on $\wt{\C}$. Since $\wt{\C}$ is complete, $\mF_u$ is convergent to some $F(u)\in\wt{\C}$. An easy adaptation of the proof of Proposition 13 in \cite[Chapter III, Section 6]{Robertson:80} to the $\wt{\C}$-module $\G$ and its algebraic dual $L(\G,\wt{\C})$ shows that $F:u\to F(u)$ is $\wt{\C}$-linear and that $\mF\to F$ according to the weak topology $\sigma(L(\G,\wt{\C}),\G)$.
\end{proof}
\begin{definition}
\label{polar_set}
Let $(\G,\mH,\bil)$ be a pairing of $\wt{\C}$-modules and $A$ be a subset of $\G$. The \emph{polar} of $A$ is the subset $A^{\circ}$ of $\mH$ of those $v\in\mH$ such that $\vert \bil(u,v)\vert_\esp\le 1$ for all $u\in A$. Similarly we define the polar of a subset of $\mH$.
\end{definition}
Some elementary properties of polar sets are collected in the following proposition.
\begin{proposition}
\label{prop_polar}
\leavevmode
\begin{trivlist}
\item[(i)] If $A_1\subseteq A_2$ then $A_2^{\circ}\subseteq A_1^{\circ}$.
\item[(ii)] The polar set of $A\subseteq\G$ is a balanced convex subset of $\mH$ closed for $\sigma(\mH,\G)$.
\item[(iii)] For all invertible $\lambda\in\wt{\C}$, $(\lambda A)^{\circ}=\lambda^{-1}A^{\circ}$. In particular $A^\circ$ is absorbent if and only if $A$ is bounded in $(\G,\sigma(\G,\mH))$.
\item[(iv)] $(\bigcup_{i\in I}A_i)^\circ =\bigcap_{i\in I}A_i^{\circ}$.
\end{trivlist}
\end{proposition}
\begin{proof}
We omit the proof of the first and the fourth assertion since it is a simple application of Definition \ref{polar_set}.

Let $A\subseteq\G$. $A^\circ$ is balanced since for all $\lambda\in\wt{\C}$ with $\vert\lambda\vert_\esp\le 1$, if $v\in A^\circ$ then $\vert\bil(u,\lambda v)\vert_\esp=\vert\lambda \bil(u,v)\vert_\esp\le \vert\lambda\vert_\esp\vert \bil(u,v)\vert_\esp\le 1$ on $A$. For each $v_1,v_2\in A^\circ$ and $u\in A$, the estimate
\[
\vert \bil(u,v_1+v_2)\vert_\esp=\vert \bil(u,v_1)+\bil(u,v_2)\vert_\esp 
\le\max\{\vert \bil(u,v_1)\vert_\esp,\vert \bil(u,v_2)\vert_\esp\}\le 1
\]
holds, i.e. $A^\circ+A^\circ\subseteq A^\circ$. This result combined with the fact that $A^\circ$ is balanced shows that $A^\circ$ is convex. Finally $A^\circ$ is closed in $(\mH,\sigma(\mH,\G))$ since it may be written as $\cap_{u\in A}A_u$ where $A_u:=\{ v\in\mH:\ \vert \bil(u,v)\vert_\esp\le 1\}$ is closed by definition of the weak topology on $\mH$. 

Take now $\lambda$ invertible in $\wt{\C}$. The equality $(\lambda A)^{\circ}=\lambda^{-1}A^{\circ}$ is guaranteed by the $\wt{\C}$-bilinearity of $\bil$. If $A^\circ$ is absorbent then for all $v\in\mH$ there exists $a\in\R$ such that $v\in[(\eps^b)_\eps]A^\circ=([(\eps^{-b})_\eps]A)^\circ$ for all $b\le a$. As a consequence, recalling that a typical neighborhood of the origin in $(\G,\sigma(\G,\mH))$ is of the form $U:=\{u\in\G:\ \max_{i=1,...,n}\vert\bil(u,v_i)\vert_\esp\le\eta\}=\{u\in\G:\ \max_{i=1,...,n}\vert\bil([(\eps^{\log\eta})_\eps]u,v_i)\vert_\esp\le 1\}$ we find $a$ in $\R$ such that $[(\eps^{-b})_\eps]A\subseteq U$ provided $b\le a+\log\eta$. This inclusion shows that $A$ is bounded for $\sigma(\G,\mH)$. Conversely if $A$ is bounded then for all $v\in\mH$ there exists $a\in\R$ such that $A$ is contained in $[(\eps^b)_\eps]\{u\in\G:\ \vert \bil(u,v)\vert_\esp\le 1\}$ for all $b\le a$. By the first statement of this proposition we conclude that $A^\circ$ absorbs every $v$ in $\mH$ since
\begin{multline*}
[(\eps^{-b})_\eps] v\in [(\eps^{-b})_\eps]\{u\in\G:\ \vert \bil(u,v)\vert_\esp\le 1\}^{\circ}\\
=([(\eps^b)_\eps]\{u\in\G:\ \vert \bil(u,v)\vert_\esp\le 1\})^\circ\subseteq A^\circ
\end{multline*}
for any $b$ smaller than $a$.
\end{proof} 
By Proposition \ref{prop_polar} the polar set of a $\sigma(\G,\mH)$-bounded subset of $\G$ is absorbent and absolutely convex. Hence its gauge defines an ultra-pseudo-seminorm on $\mH$.
\begin{definition}
\label{def_polar_topology}
Let $(\G,\mH,\bil)$ be a pairing of $\wt{\C}$-modules. A topology on $\mH$ is said to be \emph{polar} if it is determined by the family of ultra-pseudo-seminorms $\{\mP_{A^\circ}\}_{A\in\mathcal{A}}$ where $\mathcal{A}$ is a collection of $\sigma(\G,\mH)$-bounded subsets of $\G$. When $\mathcal{A}$ is the collection of all $\sigma(\G,\mH)$-bounded subsets of $\G$ then the corresponding polar topology is called \emph{strong topology} and denoted by $\beta(\mH,\G)$.
\end{definition} 
Note that $[(\eps^{-b})_\eps]v\in A^\circ$ if and only if $\sup_{u\in A}\vert \bil(u,v)\vert_\esp\le \esp^{-b}$. It follows that $\mP_{A^\circ}(v)=\sup_{u\in A}\vert\bil(u,v)\vert_\esp$ for every $\sigma(\G,\mH)$-bounded subset $A$ of $\G$. It is clear that the strong topology $\beta(\mH,\G)$ is finer than the weak topology $\sigma(\mH,\G)$.

We now deal with a particular type of locally convex topological $\wt{\C}$-modules whose topological duals have some interesting properties as we shall see in Proposition \ref{complete_dual}.
\begin{definition}
\label{def_bornivorous}
In a topological $\wt{\C}$-module $\G$ a set $S$ is said to be \emph{bornivorous} if it absorbs every bounded subset of $\G$, that is, for all bounded subsets $A$ of $\G$ there exists $a\in\R$ such that $A\subseteq [(\eps^b)_\eps]S$ for every $b\le a$.
\end{definition}
\begin{definition}
\label{def_bornological}
A locally convex topological $\wt{\C}$-module $\G$ is \emph{bornological} if every balanced, convex and bornivorous subset of $\G$ is a neighborhood of the origin.
\end{definition}
In the sequel we discuss the main result on bornological $\wt{\C}$-modules concerning bounded $\wt{\C}$-linear maps and we give some examples.
\begin{proposition}
\label{prop_continuity_born}
Let $\G$ be a bornological locally convex topological $\wt{\C}$-module and $\mH$ be an arbitrary locally convex topological $\wt{\C}$-module. If $T$ is a $\wt{\C}$-linear bounded map from $\G$ into $\mH$ then $T$ is continuous.
\end{proposition}
\begin{proof}
Let $V$ be an absolutely convex neighborhood of $0$ in $\mH$ and $A$ be a bounded subset of $\G$. By hypothesis $T(A)$ is bounded in $\mH$, i.e. there exists $a\in\R$ such that $T(A)\subseteq[(\eps^b)_\eps]V$ for all $b\le a$. $T^{-1}(V)$ is absolutely convex in $\G$ and, as proven above, it absorbs every bounded subset of $\G$. Therefore, since $\G$ is bornological, $T^{-1}(V)$ is a neighborhood of $0$ in $\G$ and $T$ is continuous.
\end{proof}
\begin{proposition}
\label{prop_examples_born}
\leavevmode
\begin{trivlist}
\item[(i)] Every locally convex topological $\wt{\C}$-module $\G$ which has a countable base of neighborhoods of the origin is bornological.
\item[(ii)] The inductive limit $\G$ of a family of bornological locally convex topological $\wt{\C}$-modules $(\G_\gamma)_{\gamma\in\Gamma}$ is bornological.
\end{trivlist}
\end{proposition}
\begin{proof}
We easily adapt the proof of the corresponding results for locally convex topological vector spaces presented in \cite[Chapter 3, Section 7, Propositions 3,4]{Horvath:66}.\\
$(i)$ Let $\G$ be a locally convex topological $\wt{\C}$-module with a countable base of neighborhoods of the origin. We may choose a balanced base of neighborhoods of the origin $(V_n)_n$ such that $V_{n+1}\subseteq V_n$ for all $n\in\N$. Let $U$ be a balanced, convex and bornivorous subset of $\G$. We want to prove that $U$ contains some $[(\eps^n)_\eps]V_n$. Assume that $U$ does not contain any $[(\eps^n)_\eps]V_n$. This means that we find a sequence $(u_n)_n$ of points $u_n\in([(\eps^n)_\eps]V_n)\cap(\G\setminus U)$. Now by construction $[(\eps^{-n})_\eps]u_n$ converges to $0$ in $\G$ and so the set $A:=\{[(\eps^{-n})_\eps]u_n,\ n\in\N\}$ is bounded. But $U$ does not absorb $A$ because if there existed $a\in\R$ such that $A\subseteq[(\eps^a)_\eps]U$ then $u_n\in[(\eps^{n+a})_\eps]U\subseteq U$ for $n$ large enough, in contradiction to our choice of the sequence $(u_n)_n$. Thus $U$ is a neighborhood of $0$ in $\G$.\\
$(ii)$ Let $\iota_\gamma:\G_\gamma\to\G$ be the family of $\wt{\C}$-linear maps which defines the inductive limit topology on $\G$ and $U$ be an absolutely convex and bornivorous subset of $\G$. Hence $\iota_\gamma^{-1}(U)$ is absolutely convex in $\G_\gamma$ and by continuity of $\iota_\gamma$, if $A_\gamma$ is bounded in $\G_\gamma$ then $\iota_\gamma(A_\gamma)$ is bounded in $\G$. By Definition \ref{def_bornivorous} there exists $a_\gamma\in\R$ such that $\iota_\gamma(A_\gamma)\subseteq [(\eps^b)_\eps]U$ for all $b\le a_\gamma$ that is $A_\gamma\subseteq[(\eps^b)_\eps]\iota_\gamma^{-1}(U)$. We have proved that $\iota_\gamma^{-1}(U)$ is balanced, convex and bornivorous and since $\G_\gamma$ is bornological, $\iota_\gamma^{-1}(U)$ is a neighborhood of $0$ in $\G_\gamma$. This tells us that $U$ is a neighborhood of $0$ in $\G$.
\end{proof}
We conclude this section by considering the pairing formed by a topological $\wt{\C}$-module $\G$ and its topological dual $\L(\G,\wt{\C})$. We know that $\L(\G,\wt{\C})$ can be endowed with the separated topologies $\sigma(\L(\G,\wt{\C}),\G)$ and $\beta(\L(\G,\wt{\C}),\G)$. Since every ultra-pseudo-seminorm defining $\sigma(\G,\L(\G,\wt{\C}))$ is continuous for the original topology $\tau$ on $\G$, $\sigma(\G,\L(\G,\wt{\C})$ is coarser than $\tau$. In some particular cases the strong topology turns $\L(\G,\wt{\C})$ into a complete topological $\wt{\C}$-module.
\begin{proposition}
\label{complete_dual}
Let $\G$ be a bornological locally convex topological $\wt{\C}$-module. Then $\L(\G,\wt{\C})$ with the topology $\beta(\L(\G,\wt{\C}),\G)$ is complete.
\end{proposition}
\begin{proof}
We shall show that every Cauchy filter $\mF$ on $\L(\G,\wt{\C})$ is convergent. First of all by the same reasoning as in Proposition \ref{complete_algebraic_dual}, for all $u\in\G$ the filter $\mF_u$ generated by $\{X_u\}_{X\in\mF}$ where $X_u:=\{Tu:\, T\in X\}$, is a Cauchy filter on $\wt{\C}$ and it converges to some $F(u)\in\wt{\C}$. The function $F:\G\to\wt{\C}:u\to F(u)$ is $\wt{\C}$-linear. Moreover $F$ is continuous on $\G$. In fact every bounded subset $A$ of $\G$ is $\sigma(\G,\L(\G,\wt{\C}))$-bounded and by definition of a Cauchy filter and the strong topology on $\L(\G,\wt{\C})$, there exists $X\in\mF$ such that $X-X\subseteq A^\circ$. This means that for all $T,T'\in X$ and for all $u\in A$ we have $\vert T(u)-T'(u)\vert_\esp\le 1$. In particular we find a constant $c>0$ such that $\vert T(u)\vert_\esp\le c$ for all $T$ in $X$ and for all $u\in A$. On the other hand given $u\in A$, $F(u)$ adheres to $\mF_u$ and therefore there exists $T'\in X$ such that $\vert F(u)-T'(u)\vert_\esp\le c$. Thus for all $u\in A$ we may write
\[
\vert F(u)\vert_\esp\le\max\{\vert F(u)-T'(u)\vert_\esp,\vert T'(u)\vert_\esp\}\le c
\]
which yields that $F(A)$ is a bounded subset of $\wt{\C}$. Since $\G$ is bornological, the bounded $\wt{\C}$-linear map $F$ is continuous.

We complete the proof by proving that $\mF$ converges to $F$ in the strong topology $\beta(\L(\G,\wt{\C}),\G)$. For all $\sigma(\G,\L(\G,\wt{\C}))$-bounded subsets $A$ of $\G$ and for all $\eta>0$ there exists $X\in\F$ such that $\vert T(u)-T'(u)\vert_\esp < \eta$ for all $u$ in $A$ and $T,T'$ in $X$. Since $\F_u\to F(u)$, for all $u\in A$ there exists $T'\in X$ such that $\vert F(u)-T'(u)\vert_\esp < \eta$. Then for all $T\in X$ and $u\in A$
\[
\vert F(u)-T(u)\vert_\esp \le\max\{\vert F(u)-T'(u)\vert_\esp , \vert T'(u)-T(u)\vert_\esp\} < \eta
\]
or in other words $\mP_{A^\circ}(F-T) < \eta$. This implies the inclusion $X\subseteq F+[(\eps^{-\log\eta})_\eps]A^\circ$. A typical neighborhood of the origin in $\beta(\L(\G,\wt{\C}),\G)$ is given by $\{T:\ \max_{i=1,...,N}\mP_{A_i^\circ}(T) < \eta\}$. Hence, from the previous considerations there exists $X=\cap_{i=1,...N}X_i$, $X_i\in\mF$, such that $X\subseteq F+[(\eps^{-\log\eta})_\eps]\cap_{i=1,...,N}A_i^\circ$, which proves our assertion.
\end{proof}
\begin{remark}
\label{remark_strong_topology}
When $\G$ is a topological $\wt{\C}$-module we can restrict the family of $\sigma(\G,\L(\G,\wt{\C}))$-bounded subsets which defines the strong topology $\beta(\L(\G,\wt{\C}),\G)$ to the family of bounded subsets of $\G$. The corresponding polar topology is called topology of uniform convergence on bounded subsets of $\G$ and denoted by $\beta_b(\L(\G,\wt{\C}),\G)$ here. Clearly $\beta_b(\L(\G,\wt{\C}),\G)$ is separated and coarser then $\beta(\L(\G,\wt{\C}),\G)$. A careful inspection of the proof of Proposition \ref{complete_dual} shows that when $\G$ is a bornological locally convex topological $\wt{\C}$-module then $\L(\G,\wt{\C})$ is complete for the topology $\beta_b(\L(\G,\wt{\C}),\G)$. In Section 3.3 we shall prove that the topologies $\beta_b(\L(\G,\wt{\C}),\G)$ and $\beta(\L(\G,\wt{\C}),\G)$ coincide for a certain particular class of ultra-pseudo-normed $\wt{\C}$-modules. The general issue concerning locally convex topological $\wt{\C}$-modules remains open.
\end{remark}
We conclude our investigation into the properties of the topological dual of a locally convex topological $\wt{\C}$-module by looking at convergent sequences. More precisely we will prove that under suitable hypotheses on $\G$, if a sequence of $\wt{\C}$-linear continuous maps $T_n:\G\to\wt{\C}$ is pointwise convergent to some $T:\G\to\wt{\C}$ then $T$ is itself an element of the dual $\L(\G,\wt{\C})$. This requires some preliminary notions concerning barrels and barrelled $\wt{\C}$-modules.  
\begin{definition}
\label{def_barreled}
Let $\G$ be a locally convex topological $\wt{\C}$-module. An absorbent, balanced, convex and closed subset of $\G$ is said to be a \emph{barrel}. A locally convex topological $\wt{\C}$-module is \emph{barrelled} if every barrel is a neighborhood of the origin.
\end{definition}
We recall that a subset $A$ of $\L(\G,\wt{\C})$ ($\G$ topological $\wt{\C}$-module) is \emph{equicontinuous at $u_0\in\G$} if for every $W$ neighborhood of the origin in $\wt{\C}$ there exists a neighborhood $U$ of $u_0$ in $\G$ such that $T(u)-T(u_0)\in W$ for all $u\in U$ and $T\in A$. $A$ is equicontinuous if it is equicontinuous at every point of $\G$.

There exists a relationship among barrels of $\G$, $\sigma(\L(\G,\wt{\C}),\G)$-bounded subsets and equicontinuous subsets of $\L(\G,\wt{\C})$.
\begin{proposition}
\label{prop_barrel}
\leavevmode
\begin{trivlist}
\item[(i)] Let $\G$ be a locally convex topological $\wt{\C}$-module. If the subset $A\subseteq\L(\G,\wt{\C})$ is $\sigma(\L(\G,\wt{\C}),\G)$-bounded then there exists a barrel $B$ in $\G$ such that $A\subseteq B^\circ$.
\item[(ii)] If $\G$ is a barrelled locally convex topological $\wt{\C}$-module then every $A\subseteq\L(\G,\wt{\C})$ which is bounded for $\sigma(\L(\G,\wt{\C}),\G)$ is equicontinuous.
\end{trivlist}
\end{proposition}
\begin{proof}
$(i)$ If $A$ is $\sigma(\L(\G,\wt{\C}),\G)$-bounded then by Proposition \ref{prop_polar} its polar $A^\circ$ is absorbent, balanced and convex. Moreover, $A^\circ =\cap_{T\in A}A_T$ where $A_T:=\{u\in\G: \vert T(u)\vert_\esp\le 1\}$ is closed in $\G$ by continuity of $T$. Hence $B:=A^\circ$ is a barrel of $\G$ such that $A\subseteq B^\circ$.\\
$(ii)$ By assertion $(i)$ every $\sigma(\L(\G,\wt{\C}),\G)$-bounded subset $A$ is contained in some $B^\circ$ where $B$ is a barrel of $\G$. This means that $\vert T([(\eps^{-\log\eta})_\eps]u)\vert_\esp\le \eta$ for all $T\in A$, $u\in B$ and $\eta>0$. Since $\G$ is barrelled $B$ is a neighborhood of $0$ and by the estimate above $A$ is an equicontinuous subset of $\L(\G,\wt{\C})$.
\end{proof} 
We now give some examples of barrelled locally convex topological $\wt{\C}$-modules.
We recall that if a \emph{Baire space} is the union of a countable family of closed subsets $S_n$ at least one set $S_n$ has nonempty interior. Baire's Theorem says that a complete metrizable topological space is a Baire space. Hence every Fr\'echet $\wt{\C}$-module is a Baire space. 
\begin{proposition}
\label{prop_baire_barrelled}
A locally convex topological $\wt{\C}$-module $\G$ which is a Baire space is barrelled.
\end{proposition}
\begin{proof}
Let $B$ be an absorbent, balanced, convex and closed subset of $\G$. Since it is absorbent we may write $\G=\cup_{n\in\N}[(\eps^{-n})_\eps]B$, where each $[(\eps^{-n})_\eps]B$ is closed in $\G$. $\G$ is a Baire space. Hence there exists some $[(\eps^{-n})_\eps]B$ with nonempty interior and from the continuity of the scalar multiplication $\G\to\G:u\to[(\eps^{-n})_\eps]u$ we conclude that $\text{int}(B)\neq\emptyset$. Let $u_0\in\text{int}(B)$. We find a neighborhood $V$ of $0$ such that $u_0+V\subseteq B$ and since $B$ is balanced $-u_0$ belongs to $B$. Hence by the convexity of $B$, $V\subseteq u_0+V-u_0\subseteq B+B\subseteq B$ which yields that $B$ is a neighborhood of $0$ in $\G$.
\end{proof}
Proposition \ref{prop_baire_barrelled} allows us to say that every Fr\'echet $\wt{\C}$-module is a barrelled locally convex topological $\wt{\C}$-module. The same conclusion holds when we consider an inductive limit procedure.
\begin{proposition}
\label{prop_barrelled_inductive}
The inductive limit $\G$ of a family of barrelled locally convex topological $\wt{\C}$-modules $(\G_\gamma)_{\gamma\in\Gamma}$ is barrelled.
\end{proposition}
The proof of Proposition \ref{prop_barrelled_inductive} is left to the reader since it is an elementary application of the definition of barrelled locally convex topological $\wt{\C}$-module and the continuity properties. Before dealing with sequences $(T_n)_n$ in the dual of a barrelled locally convex topological $\wt{\C}$-module which are pointwise convergent to some map $T:\G\to\wt{\C}$ we state a preparatory lemma.
\begin{lemma}
\label{lemma_equicontinuous}
Let $\G$ be a topological $\wt{\C}$-module, $M$ an equicontinuous subset of $\L(\G,\wt{\C})$ and $\mF$ a filter on $M$. Assume that for all $u\in \G$, the filter $\mF_u$ converges to some $F(u)\in\wt{\C}$. Then the map $F:u\to F(u)$ belongs to $\L(\G,\wt{\C})$.
\end{lemma}
\begin{proof}
As already observed in the proof of Proposition \ref{complete_algebraic_dual}, $\mF_u$, the filter generated by $\{X_u\}_{X\in\mF}$, $X_u:=\{T(u), T\in X\}$, which is convergent to $F(u)$, provides a $\wt{\C}$-linear map $F:\G\to\wt{\C}$. Since $M$ is equicontinuous we have that for all $\eta>0$ there exists a neighborhood $U$ of the origin in $\G$ such that $\vert T(u)\vert_\esp\le\eta$ for all $u\in U$ and $T\in M$. By $\mF_u\to F(u)$ it follows that for all $u\in U$ and $\eta>0$ there exists $T'\in M$ such that $\vert T'(u)-F(u)\vert_\esp\le\eta$. The estimate $\vert F(u)\vert_\esp\le\max\{ \vert T'(u)-F(u)\vert_\esp, \vert T'(u)\vert_\esp\}\le\eta$ valid on $U$ entails that $F$ is continuous at $0$ and therefore continuous on $\G$.
\end{proof} 
\begin{proposition}
\label{prop_barrelled_sequence}
Let $\G$ be a barrelled locally convex topological $\wt{\C}$-module. Let $\mF$ be a filter on $\L(\G,\wt{\C})$ which contains a $\sigma(\L(\G,\wt{\C}),\G)$-bounded subset of $\L(\G,\wt{\C})$. Assume that the filter $\mF_u$ converges to some $F(u)\in\wt{\C}$. Then the map $F:u\to F(u)$ belongs to $\L(\G,\wt{\C})$.
\end{proposition}  
\begin{proof}
By assertion $(ii)$ in Proposition \ref{prop_barrel} we know that $\mF$ contains an equicontinuous subset $M$ of $\L(\G,\wt{\C})$. Let $\mO:=\{Y\subseteq M:\, \exists X\in\mF\ X\cap M\subseteq Y\}$ be the filter induced by $\mF$ on $M$. $\mO_u$ converges to $F(u)$ for all $u\in\G$. An application of Lemma \ref{lemma_equicontinuous} proves that $F$ is a $\wt{\C}$-linear and continuous map on $\G$.
\end{proof}
\begin{corollary}
\label{cor_barrelled_sequence}
Let $\G$ be a barrelled locally convex topological $\wt{\C}$-module. Suppose that $(T_n)_n$ is a sequence in $\L(\G,\wt{\C})$ such that for every $u\in\G$ the sequence $(T_n(u))_n$ converges to some $T(u)\in\wt{\C}$. Then $T:\G\to\wt{\C}:u\to T(u)$ belongs to $\L(\G,\wt{\C})$.
\end{corollary}
\begin{proof}
The elementary filter associated with the sequence $(T_n)_n$ i.e. the filter $\mF$ generated by $X_N:=\{T_n,\ n\ge N\}$, $N\in\N$, contains a $\sigma(\L(\G,\wt{\C}),\G)$-bounded subset of $\L(\G,\wt{\C})$. In fact by $T_n(u)\to T(u)$ each $X_N$ is $\sigma(\L(\G,\wt{\C}),\G)$-bounded and by construction $\mF_u\to T(u)$ for all $u\in\G$. By Proposition \ref{prop_barrelled_sequence} we conclude that $T\in\L(\G,\wt{\C})$.
\end{proof}
\begin{corollary}
\label{cor_quasi_complete}
If $\G$ is a barrelled locally convex topological $\wt{\C}$-module then the topological dual $\L(\G,\wt{\C})$ endowed with the weak-topology $\sigma(\L(\G,\wt{\C}),\G)$ is quasi-complete.
\end{corollary}
\begin{proof}
We have to show that every closed and bounded subset $M$ of $\L(\G,\wt{\C})$ is complete for the topology $\sigma(\L(\G,\wt{\C}),\G)$. Let $\mF$ be a Cauchy filter on $M$. $\mF$ generates a Cauchy filter $\mF'$ on $\L(\G,\wt{\C})$ which contains a  $\sigma(\L(\G,\wt{\C}),\G)$-bounded subset of $\L(\G,\wt{\C})$ and a Cauchy filter $\mF''$ on $(L(\G,\wt{\C}), \sigma(L(\G,\wt{\C}),\G))$. By Proposition \ref{complete_algebraic_dual} there exists $F\in L(\G,\wt{\C})$ such that $\mF''\to F$ in $L(\G,\wt{\C})$ and consequently $\mF'_u\to F(u)$ for all $u\in\G$. At this point Proposition \ref{prop_barrelled_sequence} allows to conclude that $F\in\L(\G,\wt{\C})$ and $\mF'\to F$ in $\L(\G,\wt{\C})$. Since $\mF$ is a filter on $M$ and $M$ is $\sigma(\L(\G,\wt{\C}),\G)$-closed, $F$ itself belongs to $M$ and $\mF\to F$.  
\end{proof}  
\section{$\wt{\C}$-modules of generalized functions based on a locally convex topological vector space}
In this part of the paper we focus our attention on a relevant class of examples of $\wt{\C}$-modules, whose general theory was developed in the previous sections. In the literature there already exist papers on topologies, generalized functions and applications cf. \cite{Biagioni:90, BiagCol:86, O:91} which consider spaces of generalized functions $\G_E$ based on a locally convex topological vector space $E$ and define topologies in terms of valuations and ultra-pseudo-seminorms \cite{NPS:98, Scarpalezos:92, Scarpalezos:98, Scarpalezos:00, TodorovW:00}. The topological background provided by Sections 1 and 2 allows us to consider $\G_E$ as an element of the larger family of locally convex topological $\wt{\C}$-modules and to deal with issues as boundedness, completeness and topological duals. For the sake of exposition we organize the following notions and results in three subsections.
\subsection{Definition and basic properties of $\G_E$}
\begin{definition}
\label{defGE}
Let $E$ be a locally convex topological vector space topologized through the family of seminorms $\{p_i\}_{i\in I}$. The elements of  
\beq
\label{defME}
 \M_E := \{(u_\eps)_\eps\in E^{(0,1]}:\, \forall i\in I\,\, \exists N\in\N\quad p_i(u_\eps)=O(\eps^{-N})\, \text{as}\, \eps\to 0\}
\eeq
and
\beq
\label{defNE}
 \Neg_E := \{(u_\eps)_\eps\in E^{(0,1]}:\, \forall i\in I\,\, \forall q\in\N\quad p_i(u_\eps)=O(\eps^{q})\, \text{as}\, \eps\to 0\},
\eeq  
are called $E$-moderate and $E$-negligible, respectively. We define the space of \emph{generalized functions based on $E$} as the factor space $\G_E := \M_E / \Neg_E$.
\end{definition}
It is clear that the definition of $\G_E$ does not depend on the family of seminorms which determines the locally convex topology of $E$. We adopt the notation $u=[(u_\eps)_\eps]$ for the class $u$ of $(u_\eps)_\eps$ in $\G_E$ and we embed $E$ into $\G_E$ via the constant embedding $f\to[(f)_\eps]$. By the properties of seminorms on $E$ we may define the product between complex generalized numbers and elements of $\G_E$ via the map  $\wt{\C}\times\G_E\to\G_E:([(\lambda_\eps)_\eps],[(u_\eps)_\eps])\to[(\lambda_\eps u_\eps)_\eps]$, which equips $\G_E$ with the structure of a $\wt{\C}$-module.

Since the growth in $\eps$ of an $E$-moderate net is estimated in terms of any seminorm $p_i$ of $E$, it is natural to introduce the \emph{$p_i$-valuation} of $(u_\eps)_\eps\in\M_E$ as 
\beq
\label{valpi}
\val_{p_i}((u_\eps)_\eps) := \sup\{ b\in\R:\quad p_i(u_\eps)=O(\eps^b)\quad \text{as}\, \eps\to 0\}.
\eeq
Note that $\val_{p_i}((u_\eps)_\eps)=\val((p_i(u_\eps))_\eps)$ where the function $\val$ in \eqref{val_M} gives the valuation on $\wt{\C}$. Clearly $\val_{p_i}$ maps $\M_E$ into $(-\infty,+\infty]$ and the following properties hold:
\begin{trivlist}
\item[(i)] $\val_{p_i}((u_\eps)_\eps)= +\infty$ for all $i\in I$ if and only if $(u_\eps)_\eps\in\Neg_E$,
\item[(ii)] $\val_{p_i}((\lambda_\eps u_\eps)_\eps)\ge \val((\lambda_\eps)_\eps)+\val_{p_{i}}((u_\eps)_\eps)$ for all $(\lambda_\eps)_\eps\in\E_M$ and $(u_\eps)_\eps\in\M_E$,
\item[(ii)'] $\val_{p_i}((\lambda_\eps u_\eps)_\eps)= \val((\lambda_\eps)_\eps)+\val_{p_{i}}((u_\eps)_\eps)$ for all $(\lambda_\eps)_\eps=(c\eps^b)_\eps$, $c\in\C$, $b\in\R$,
\item[(iii)] $\val_{p_i}((u_\eps)_\eps + (v_\eps)_\eps)\ge \min\{ \val_{p_i}((u_\eps)_\eps), \val_{p_i}((v_\eps)_\eps)\}$.
\end{trivlist}
Assertion $(i)$ combined with $(iii)$ shows that $\val_{p_i}((u_\eps)_\eps)=\val_{p_i}((u'_\eps)_\eps)$ if $(u_\eps-u'_\eps)_\eps$ is $E$-negligible. This means that we can use \eqref{valpi} for defining the \emph{$p_i$-valuation} $\val_{p_i}(u)=\val_{p_i}((u_\eps)_\eps)$ of a generalized function $u=[(u_\eps)_\eps]\in\G_E$.

$\val_{p_i}$ is a valuation in the sense of Definition \ref{def_ultra_pseudo} and thus $\mP_i(u):=\esp^{-\val_{p_i}(u)}$ is an ultra-pseudo-seminorm on the $\wt{\C}$-module $\G_E$. By Theorem \ref{theorem_convex}, $\G_E$ endowed with the topology of the ultra-pseudo-seminorms $\{\mP_i\}_{i\in I}$ is a locally convex topological $\wt{\C}$-module. Following \cite{NPS:98, Scarpalezos:92, Scarpalezos:98, Scarpalezos:00} we use the adjective ``sharp'' for the topology induced by the ultra-pseudo-seminorms $\{\mP_i\}_{i\in I}$. The sharp topology on $\G_E$, here denoted by $\tau_\sharp$ is independent of the choice of the family of seminorms which determines the original locally convex topology on $E$. The structure of the subspace $\Neg_E$ has some interesting influence on $\tau_\sharp$.
\begin{proposition}
\label{prop_separated_sharp}
$(\G_E,\tau_\sharp)$ is a separated locally convex topological $\wt{\C}$-module.
\end{proposition}
\begin{proof}
By definition of $\Neg_E$ if $u\neq 0$ in $\G_E$ then $\val_{p_i}((u_\eps)_\eps)\neq +\infty$ for some $i\in I$. This means that $\mP_i(u)>0$ and by Proposition \ref{prop_T2_convex} $\tau_\sharp$ is a separated topology.
\end{proof}
\begin{proposition}
\label{prop_ideal_topology}
Let $E$ be a locally convex topological vector space.
\leavevmode
\begin{trivlist}
\item[(i)] If $E$ is topologized through an increasing sequence $\{p_i\}_{i\in\N}$ of seminorms and 
\beq
\label{charac_NE}
\Neg_E =\{ (u_\eps)_\eps \in\M_E:\, \forall q\in\N\quad p_0(u_\eps)=O(\eps^q)\,\, \text{as}\, \eps\to 0\},
\eeq
then each $\mP_i$ is an ultra-pseudo-norm on $\G_E$.
\item[(ii)] If $E$ has a countable base of neighborhoods of the origin then $\G_E$ with the sharp topology is metrizable.
\end{trivlist}
\end{proposition}
\begin{proof}
Concerning the first assertion we have to prove that $\mP_i(u)=0$ implies $u=0$ in $\G_E$. From $\mP_i(u)=0$ it follows that $p_i(u_\eps)=O(\eps^q)$ for all $q\in\N$. Since $p_0(u_\eps)\le p_i(u_\eps)$, \eqref{charac_NE} leads to $(u_\eps)_\eps\in\Neg_E$. Combining Proposition \ref{prop_separated_sharp} with the assumption $(ii)$, we obtain that $\G_E$ is a separated locally convex topological $\wt{\C}$-module with a countable base of neighborhoods of the origin. Hence by Theorem \ref{theorem_metric} it is metrizable.
\end{proof}
\begin{proposition}
\label{complete_G_E}
If $E$ is a locally convex topological vector space with a countable base of neighborhoods of the origin then $\G_E$ with the sharp topology $\tau_\sharp$ is complete.
\end{proposition}
This result, in terms of convergence of Cauchy sequences, was already proven in \cite[Proposition 2.1]{Scarpalezos:98}. For the convenience of the reader we add some details to the sketch of the proof given there.
\begin{proof}
As shown in Proposition \ref{prop_ideal_topology} $\G_E$ is metrizable and therefore by Remark \ref{remark_complete} it is sufficient to prove that any Cauchy sequence in $\G_E$ is convergent. It is not restrictive to assume that $E$ is topologized through an increasing sequence of seminorms $\{p_k\}_{k\in\N}$. If $(u_n)_n$ is a Cauchy sequence in $\G_E$ we may extract a subsequence $(u_{n_k})_k$ such that $\val_{p_k}((u_{n_{k+1},\eps}-u_{n_k,\eps})_\eps)>k$ for all $k\in\N$. As in the proof of Proposition \ref{prop_C_complete} we obtain a decreasing sequence $\eps_k\searrow 0$, $\eps_k\le 2^{-k}$ such that $p_k(u_{n_{k+1},\eps}-u_{n_k,\eps})\le\eps^k$ for all $\eps\in(0,\eps_k)$. Let   
\[
h_{k,\eps}=\begin{cases} u_{n_{k+1},\eps}-u_{n_k,\eps}\  & \eps\in(0,\eps_k),\\
0\ & \eps\in[\eps_k,1].
\end{cases}
\] 
Obviously $(h_{k,\eps})_\eps$ belongs to $\M_E$ and for all $k'\le k$, $p_{k'}(h_{k,\eps})\le\eps^k$ on the interval $(0,1]$. The sum $u_\eps:=u_{n_0,\eps}+	\sum_{k=0}^\infty h_{k,\eps}$ is locally finite and $E$-moderate since for all $\overline{k}\in\N$
\begin{multline*}
p_{\overline{k}}(u_\eps)\le p_{\overline{k}}(u_{n_0,\eps})+\sum_{k=0}^{\overline{k}}p_{\overline{k}}(h_{k,\eps})+\sum_{k=\overline{k}+1}^\infty p_{\overline{k}}(h_{k,\eps})\\
\le p_{\overline{k}}(u_{n_0,\eps})+\sum_{k=0}^{\overline{k}}p_{\overline{k}}(h_{k,\eps})+\sum_{k=\overline{k}+1}^\infty \frac{1}{2^k}.
\end{multline*}
Finally for all $\overline{k}\ge 1$ and for all $\eps\in(0,\eps_{\overline{k}-1})$
\beq
\label{estimate_pk}
p_{\overline{k}}(u_{n_{\overline{k}},\eps}-u_\eps)= p_{\overline{k}}\big(-\sum_{k=\overline{k}}^\infty h_{k,\eps}\big)\le \sum_{k=\overline{k}}^\infty \eps^{\overline{k}-1}\eps_k\le \eps^{\overline{k}-1}\sum_{k=\overline{k}}^\infty\frac{1}{2^k}.
\eeq
By \eqref{estimate_pk} we conclude that for all $\overline{k}\ge 1$, for all $q\in\N$, for all $k\ge\max\{\overline{k},q+1\}$ there exists $\eta\in(0,1]$ such that $p_{\overline{k}}(u_{n_k,\eps}-u_\eps)\le\eps^q$ on $(0,\eta]$. In other words $(u_{n_k})_k$ is convergent to $u$ in $\G_E$. Consequently $(u_n)_n$ itself converges to $u$ in $\G_E$.
\end{proof}
\begin{remark}
\label{remark_topological_algebras}
We recall that a 
$\wt{\C}$-module $\G$ is a \emph{$\wt{\C}$-algebra} if there is given a multiplication $\G\times\G\to\G:(u,v)\to uv$ such that $(uv)w=u(vw)$, $u(v+w)=uv+uw$, $(u+v)w=uw+vw$, $\lambda(uv)=(\lambda u)v=u(\lambda v)$ for all $u,v,w$ in $\G$ and $\lambda\in\wt{\C}$. In analogy with the theory of topological algebras we say that a $\wt{\C}$-algebra $\G$ is a topological $\wt{\C}$-algebra if it is equipped with a $\wt{\C}$-linear topology which makes the multiplication on $\G$ continuous. As an explanatory example let us consider an algebra $E$ on $\C$ and a family of seminorms $\{p_i\}_{i\in I}$ on $E$. Assume that for all $i\in I$ there exist finite subsets $I_0,I_0'$ of $I$ and a constant $C_i>0$ such that for all $u,v\in E$
\beq
\label{estimate_seminorm_product}
p_i(uv)\le C_i\max_{j\in I_0}p_j(u)\ \max_{j\in I'_0}p_j(v).
\eeq
Then $\G_E$ with the sharp topology determined by the ultra-pseudo-seminorms $\{\mP_i\}_{i\in I}$ is a locally convex topological $\wt{\C}$-module and a topological $\wt{\C}$-algebra since from \eqref{estimate_seminorm_product} it follows that $$\mP_i(uv)\le \max_{j\in I_0}\mP_j(u)\ \max_{j\in I'_0}\mP_j(v)$$ for all $i\in I$.
\end{remark}
In the sequel we collect some examples of locally convex topological $\wt{\C}$-modules which occur in Colombeau theory. For details and explanations about Co\-lom\-be\-au generalized functions we mainly refer to \cite{Colombeau:85, GKOS:01, NPS:98}.
\begin{example}
\label{example_GE}
\bf{Colombeau algebras obtained as $\wt{\C}$-modules $\G_E$}\rm\\
Particular choices of $E$ in Definition \ref{defGE} give us known algebras of generalized functions and the corresponding sharp topologies. This is of course the case for $E=\C$ and $\G_E=\wt{\C}$ which is an ultra-pseudo-normed $\wt{\C}$-module and more precisely a topological $\wt{\C}$-algebra. 

Consider now an open subset $\Om$ of $\R^n$. $E=\E(\Om)$, i.e. the space $\Cinf(\Om)$ topologized through the family of seminorms $p_{K_i,j}(f)=\sup_{x\in K_i, |\alpha|\le j}|\partial^\alpha f(x)|$, where $K_0\subset K_1\subset.... K_i\subset...$ is a countable and exhausting sequence of compact subsets of $\Om$, provides $\G_E=\G(\Om)$ (\cite{Colombeau:85,GKOS:01}). By Propositions \ref{prop_separated_sharp}, \ref{complete_G_E} and Remark \ref{remark_topological_algebras}, $\G(\Om)$ endowed with the sharp topology determined by $\{\mP_{K_i,j}\}_{i\in\N, j\in\N}$ is a Fr\'echet $\wt{\C}$-module and a topological $\wt{\C}$-algebra. Other examples of Fr\'echet $\wt{\C}$-modules which are also topological $\wt{\C}$-algebras are given by $\G_E$ when $E$ is $\S(\R^n)$ or $W^{\infty,p}(\R^n)$, $p\in[1,+\infty]$.
In this way we construct the algebras $\GS(\R^n)$ (\cite[Definition 2.10]{GGO:03}) and $\G_{p,p}(\R^n)$ (\cite{BO:92}) respectively, whose sharp topologies are obtained from $p_k(f)=\sup_{x\in\R^n, |\alpha|\le k}(1+|x|)^k|\partial^\alpha f(x)|$, $f\in\S(\R^n)$ and $q_k(g)=\max_{|\alpha|\le k}\Vert\partial^\alpha g\Vert_p$, $g\in W^{\infty,p}(\R^n)$, with $k$ varying in $\N$. In \cite{Garetto:04c} we prove that a characterization as \eqref{charac_NE} holds for the ideals $\Neg_{\S(\R^n)}=\NS(\R^n)$ and $\Neg_{W^{\infty,p}(\R^n)}=\Neg_{p,p}(\R^n)$. As a consequence $\mP_k$ and $\mQ_k$ are ultra-pseudo-norms on $\GS(\R^n)$ and $\G_{p,p}(\R^n)$ respectively.
\end{example}

We concentrate now on the subalgebra $\G_c(\Om)$ of generalized functions in $\G(\Om)$ with compact support. It will turn out that $\G_c(\Om)$ can be equipped with a strict inductive limit topology, but this procedure requires some preliminary investigations. For technical reason we begin by recalling the basic notions of point value theory in Colombeau algebras \cite{GKOS:01, OK:99}, which will be used in the sequel. 

The set of generalized points $\wt{x}\in\wt{\Om}$ is defined as the factor $\Om_M /\sim$, where $\Om_M:=\{(x_\eps)_\eps\in\Om^{(0,1]}:\, \exists N\in\N\ |x_\eps|=O(\eps^{-N})\, \text{as}\, \eps\to 0\}$ and $\sim$ is the equivalence relation given by 
\[
(x_\eps)_\eps\sim(y_\eps)_\eps\ \Leftrightarrow\ \forall q\in\N\quad |x_\eps-y_\eps|=O(\eps^q).
\]
We say that $\wt{x}\in\wt{\Om}_{\rm{c}}$ if it has a representative $(x_\eps)_\eps$ such that $x_\eps$ belongs to a compact set $K$ of $\Om$ for small $\eps$. One can show that the generalized point value of $u\in\G(\Om)$ at $\wt{x}\in\wt{\Om}_{\rm{c}}$,
\[
u(\wt{x}):=[(u_\eps(x_\eps))_\eps]
\]
is a well-defined element of $\wt{\C}$ and Theorem 1.2.46 in \cite{GKOS:01} allows the following characterizations of generalized functions in terms of their point values:
\beq
\label{point_value_charac}
u=0\ \text{in}\ \G(\Om)\qquad \Leftrightarrow\qquad \forall \wt{x}\in\wt{\Om}_{\rm{c}}\quad u(\wt{x})=0\ \text{in}\ \wt{\C}.
\eeq
\begin{example}
\label{example_Gc}
\bf{The Colombeau algebra of compactly supported generalized functions}\rm\\
For $K\Subset\Omega$ we denote by $\G_K(\Om)$ the space of all generalized functions in $\G(\Om)$ with support contained in $K$. Note that $\G_K(\Om)$ is contained in $\G_{\D_{K'}(\Om)}$ for all compact subsets $K'$ of $\Om$ such that $K\subseteq {\rm{int}}(K')$, where $\D_{K'}(\Om)$ is the space of all smooth functions $f$ with $\supp\, f\subseteq K'$. In fact if $\text{supp}\, u\subseteq K$ we can always find a $\D_{K'}(\Om)$-moderate representative $(u_\eps)_\eps$ and for all representatives $(u_\eps)_\eps, (u'_\eps)_\eps$ of this type, $(u_\eps -u'_\eps)_\eps\in\Neg_{\D_{K'}(\Om)}$. With this choice of $(u_\eps)_\eps$, from $\Neg_{\D_{K'}(\Om)}\subseteq \Neg(\Om)$ we have that
\[
\G_K(\Om)\to \G_{\D_{K'}(\Om)}:u\to (u_\eps)_\eps +\Neg_{\D_{K'}(\Om)}
\]
is a well-defined and injective $\wt{\C}$-linear map. Moreover, by $\M_{\D_{K'}(\Om)}\cap\Neg(\Om)\subseteq\Neg_{\D_{K'}(\Om)}$, $\G_{\D_{K'}(\Om)}$ is naturally embedded into $\G(\Om)$ via $$\G_{\D_{K'}(\Om)}\to\G(\Om):(u_\eps)_\eps+\Neg_{\D_{K'}(\Om)}\to (u_\eps)_\eps +\Neg(\Om).$$
In $\G(\Om)$ the $p_{K,n}$-valuation where $p_{K,n}(f)=\sup_{x\in K, |\alpha|\le n}|\partial^\alpha f(x)|$ is obtained as the valuation of the complex generalized number $\sup_{x\in K, |\alpha|\le n}|\partial^\alpha u(x)|:=(\sup_{x\in K, |\alpha|\le n}|\partial^\alpha u_\eps(x)|)_\eps +\Neg$. Hence for $K,K'\Subset\Om$, $K\subseteq\text{int}(K')$,
\beq
\label{valGK}
\val_{K,n}(u)=\val_{p_{K',n}}(u)
\eeq
is a valuation on $\G_K(\Om)$. More precisely \eqref{valGK} does not depend on $K'$ since for any $K'_1,K'_2$ containing $K$ in their interiors and for any $u\in\G_K(\Om)$ we have that 
\[
\val_{p_{K'_1,n}}(u)\ge\inf\big\{\val_{p_{K_1'\setminus\text{int}(K'_1\cap K'_2),n}}(u),\val_{p_{K'_2,n}}(u)\big\}
= \val_{p_{K'_2,n}}(u).
\]
$\G_K(\Om)$ with the topology induced by the ultra-pseudo-seminorms $\{\mP_{\G_K(\Om),n}(u)\\:=\esp^{-\val_{K,n}(u)}\}_{n\in\N}$ is a locally convex topological $\wt{\C}$-module and by construction its topology coincides with the topology induced by any $\G_{\D_{K'}(\Om)}$ with $K\subseteq \text{int}(K')$. In particular, the $\wt{\C}$-module $\G_K(\Om)$ is separated and by Theorem \ref{theorem_metric} it is metrizable. Finally assume that $u\in\G_{\D_{K'}(\Om)}$ adheres to $\G_K(\Om)$. We find a sequence $(u_n)_n\in\G_K(\Om)$ such that $\val_{p_{K',0}}(u-u_n)\ge n$ for all $n\in\N$. Recall that for all $\wt{x}\in\wt{V}_{\rm{c}}$, where $V=\Omega\setminus K$, the point values $u_n(\wt{x})$ are zero in $\wt{\C}$ and 
\[
\val_{\wt{\C}}(u(\wt{x})) \ge \min\{\val_{p_{K',0}}(u-u_n),\val_{\wt{\C}}(u_n(\wt{x}))\} = \val_{p_{K',0}}(u-u_n).
\]
Consequently $u(\wt{x})=0$ in $\wt{\C}$ and by \eqref{point_value_charac} $\text{supp}\, u\subseteq K$. We just proved that $\G_K(\Om)$ is closed in $\G_{\D_{K'}(\Om)}$ and since $\G_{\D_{K'}(\Om)}$ with its sharp topology is complete, Remark \ref{remark_complete} allows to conclude that \emph{$(\G_{K}(\Om),\{\mP_{\G_K(\Om),n}\}_{n\in\N})$ is a Fr\'echet $\wt{\C}$-module}.\\
Note that if $K_1\subseteq K_2$ then $\G_{K_1}(\Om)\subseteq\G_{K_2}(\Om)$ and that $\G_{K_2}(\Om)$ induces on $\G_{K_1}(\Om)$ the original topology. By construction 
$\G_{K_1}(\Om)$ is closed in $\G_{K_2}(\Om)$.

Let $(K_n)_{n\in\N}$ be an exhausting sequence of compact subsets of $\Omega$ such that $K_n\subseteq K_{n+1}$. Clearly $\G_c(\Om)=\cup_{n\in\N}\G_{K_n}(\Om)$. Each $\G_{K_n}(\Om)$ is a Fr\'echet $\wt{\C}$-module and the assumptions of Definition \ref{def_strict}, Theorem \ref{theorem_bounded_strict} and Theorem \ref{theorem_strict_complete} are satisfied by $\G_n=\G_{K_n}(\Om)$. Therefore $\G_c(\Om)$ endowed with the strict inductive limit topology of the sequence $(\G_{K_n}(\Om))_{n}$ is a separated and complete locally convex topological $\wt{\C}$-module. Obviously this topology is independent of the choice of the covering $(K_n)_n$.\\
Applying Corollary \ref{cor_cont_strict} to this context we have that \emph{a sequence $(u_n)_n$ of generalized functions with compact support converges to $0$ in $\G_c(\Om)$ if and only if it is contained in some $\G_K(\Om)$ and convergent to $0$ there.}
\end{example}
\begin{remark}
Classically \cite[Example 7, p.170]{Horvath:66} the topology on $\D(\Om)$ is determined by the seminorms
\[
p_\theta(f)=\sup_{\alpha\in\N^n}\sup_{x\in\Om}|\theta_\alpha(x)\partial^\alpha f(x)|,
\] 
where $\theta$ runs through all possible families $\theta=(\theta_\alpha)_{\alpha\in\N^n}$ of continuous functions on $\Om$ with $({\rm{supp}}\,\theta_\alpha)_{\alpha\in\N^n}$ locally finite. Easy computations show that 
\beq
\label{def_iotaD}
\iota_{\D}:\Gc(\Om)\to\G_{\D(\Om)}:u\to (u_\eps)_\eps +\Neg_{\D(\Om)},
\eeq
where $(u_\eps)_\eps$ is any representative of $u$ with supp$(u_\eps)$ contained in the same compact set of $\Om$ for all $\eps\in(0,1]$, is well-defined and injective.
One may think of endowing $\G_c(\Om)$ with the locally convex $\wt{\C}$-linear topology induced by the sharp topology on $\G_{\D(\Om)}$ via $\iota_\D$. Denoting this topology by $\tau_\D$ and the strict inductive limit topology of Example \ref{example_Gc} by $\tau$, we have that $\tau_\D$ is coarser than $\tau$ since every embedding $\G_{K_n}(\Om)\to\G_c(\Om)$ is continuous for $\tau_\D$ on $\G_c(\Om)$. In detail, taking $K'_n\Subset\Om$ with $K_n\subseteq\text{int}(K'_n)$ for all $(\theta_\alpha)_\alpha$ there exists $N\in\N$ and $C>0$ such that the estimate
\beq
\label{est_ptheta_cont}
\begin{split}
\sup_{\alpha\in\N^n}\sup_{x\in\Om}|\theta_\alpha(x)\partial^\alpha u_\eps(x)|&=\sup_{|\alpha|\le N}\sup_{x\in K'_n}|\theta_\alpha(x)\partial^\alpha u_\eps(x)|\\
&\le C\sup_{|\alpha|\le N, x\in K'_n}|\partial^\alpha u_\eps(x)|= C\, p_{K'_n,N}(u_\eps)
\end{split}
\eeq
holds for every representative of $u\in\G_{K_n}(\Om)$ with $\text{supp}(u_\eps)\subseteq K'_n$ for all $\eps\in(0,1]$. \eqref{est_ptheta_cont} implies
\[
\mP_\theta(u)\le \mP_{\G_{K_n}(\Om),N}(u),\qquad u\in\G_{K_n}(\Om)
\]
and by Corollary \ref{corollary_linear} guarantees the continuity of the embeddings mentioned above.

In general $\tau_\D$ does not coincide with $\tau$. This is shown by the fact that there exist sequences in $\G_c(\R^n)$ which converge to $0$ with respect to $\tau_\D$ but not with respect to $\tau$. Indeed, for every $\psi\in\Cinfc(\R^n)$ the generalized functions $u_n:=(\eps^n\psi(\frac{x}{n}))_\eps+\Neg(\R^n)$ have compact support and at fixed $n$
\[
\sup_{\alpha\in\N^n}\sup_{x\in\R^n}\big|\,\eps^n n^{-|\alpha|}\theta_\alpha(x)\partial^\alpha\psi(\frac{x}{n})\big|=O(\eps^n),\quad \text{as}\, \eps\to 0
\]
Thus $\val_{p_\theta}(u_n)\ge n\to +\infty$. This means that $(u_n)_n$ is $\tau_\D$-convergent to $0$. Since $\supp\, u_n= n\, \supp\, \psi$, by Example \ref{example_Gc} $(u_n)_n$ cannot be $\tau$-convergent to $0$.
\end{remark}
\begin{example}
\label{example_tempered}
\bf{The algebra of tempered generalized functions $\Gt(\R^n)$}\rm\\
The algebra of tempered generalized functions $\Gt(\R^n)$ may be introduced referring to the constructions of \cite{Colombeau:85, GKOS:01} as the factor space $\Et(\R^n)/\Nt(\R^n)$, where $\Et(\R^n)$ is the algebra of all  \emph{$\tau$-moderate} nets $(u_\eps)_\eps\in\Et[\R^n]:=\mO_M(\R^n)^{(0,1]}$ such that 
\beq
\label{tau_moderate}
\forall \alpha\in\N^n\, \exists N\in\N\qquad \sup_{x\in\R^n}(1+|x|)^{-N}|\partial^\alpha u_\eps(x)|=O(\eps^{-N})\qquad \text{as}\ \eps\to 0
\eeq
and $\Nt(\R^n)$ is the ideal of all \emph{$\tau$-negligible} nets $(u_\eps)_\eps\in\Et[\R^n]$ such that
\beq
\label{tau_negligible}
\forall \alpha\in\N^n\, \exists N\in\N\, \forall q\in\N\quad \sup_{x\in\R^n}(1+|x|)^{-N}|\partial^\alpha u_\eps(x)|=O(\eps^{q})\ \text{as}\ \eps\to 0.
\eeq
Theorem 1.2.25 in \cite{GKOS:01} shows that $\Nt(\R^n)$ coincides with the set of all $(u_\eps)_\eps\in\Et(\R^n)$ whose $0$-th derivative satisfies \eqref{tau_negligible} i.e. $\exists N\in\N\, \forall q\in\N$\ $\sup_{x\in\R^n}(1+|x|)^{-N}|u_\eps(x)|=O(\eps^{q})$. Moreover for each $\wt{x}\in\wt{\R^n}$, $u(\wt{x}):=[(u_\eps(x_\eps))_\eps]$ is a well-defined element of $\wt{\C}$ and the point value characterization  
\beq
\label{point_value_temp}
u=0\ \text{in}\ \Gt(\R^n)\qquad \Leftrightarrow\qquad \forall \wt{x}\in\wt{\R^n}\quad u(\wt{x})=0\ \text{in}\ \wt{\C}
\eeq
holds. We present a locally convex $\wt{\C}$-linear topology on $\Gt(\R^n)$ whose construction involves a countable family of different algebras of generalized functions. We denote by $\GtS(\R^n)$ the factor algebra $\Et(\R^n)/\NS(\R^n)$ (c.f. Definition 2.8 \cite{Garetto:04}) where $\NS(\R^n)=\Neg_{\S(\R^n)}$. Inspired by the definition of $\tau$-moderate nets we introduce the set
\begin{multline*}
\E^m_N(\R^n):=\{(u_\eps)_\eps\in\Et[\R^n]:\ \exists b\in\R\\ \sup_{x\in \R^n,|\alpha|\le m}(1+|x|)^{-N}|\partial^\alpha u_\eps(x)|=O(\eps^b)\ \text{as}\ \eps\to 0\}
\end{multline*}
and the $\wt{\C}$-module $\G^m_{N,\S}(\R^n):=\E^m_N(\R^n)/\NS(\R^n)$. Thus, setting $\G^m_{\tau,\S}(\R^n)\hskip-2pt:=\cup_{N\in\N}\G^m_{N,\S}(\R^n)$ we have that
\[
\GtS(\R^n)=\bigcap_{m\in\N}\G^m_{\tau,\S}(\R^n)=\bigcap_{m\in\N}\bigcup_{N\in\N}\G^m_{N,\S}(\R^n).
\]
We begin by endowing $\GtS(\R^n)$ with a locally convex $\wt{\C}$-linear topology considering $\Gt(\R^n)$ only in a second step. Every $\G^m_{N,\S}(\R^n)$ is a locally convex topological $\wt{\C}$-module for the ultra-pseudo-seminorm $\mP^m_N$ determined by the well-defined valuation 
\[
\val^m_N(u):=\sup\{b\in\R:\  \sup_{x\in \R^n,|\alpha|\le m}(1+|x|)^{-N}|\partial^\alpha u_\eps(x)|=O(\eps^b)\ \text{as}\ \eps\to 0\}.
\]
Hence, we equip $\G^m_{\tau,\S}(\R^n)$ with the inductive limit topology of the sequence $(\G^m_{N,\S}(\R^n))_{N\in\N}$ 
and we take the initial topology on $\GtS(\R^n)=\cap_{m\in\N}\G^m_{\tau,\S}(\R^n)$. 
Finally we topologize $\Gt(\R^n)$ through the finest locally convex $\wt{\C}$-linear topology such that the map
\[
\iota_{\tau,\S}:\GtS(\R^n)\to\Gt(\R^n):(u_\eps)_\eps+\NS(\R^n)\to (u_\eps)_\eps+\Nt(\R^n)
\]
is continuous. The fact that this topology is separated follows from the continuous embedding of $\Gt(\R^n)$ into the separated locally convex topological $\wt{\C}$-module $\L(\GS(\R^n),\wt{\C})$ studied in \cite{Garetto:04c}.
\end{example}
\begin{example}
\label{example_regular}
\bf{Regular generalized functions based on $E$}\rm  

For any locally convex topological vector space $(E,\{p_i\}_{i\in I})$ the set
\beq
\label{MinfE}
\M^\infty_E:=\{(u_\eps)_\eps\in E^{(0,1]}:\ \exists N\in\N\, \forall i\in I\quad p_i(u_\eps)=O(\eps^{-N})\ \text{as}\ \eps\to 0\}
\eeq
is a subspace of the set $\M_E$ of $E$-moderate nets. Therefore the corresponding factor space $\G^\infty_E:=\M^\infty_E /\Neg_E$ is a subspace of $\G_E$ whose elements are called \emph{regular generalized functions based on $E$}. When $E$ is in addition a topological algebra, i.e. estimate \eqref{estimate_seminorm_product} is satisfied, $\G_E$ and $\G^\infty_E$ are both algebras. We want to equip $\G^\infty_E$ with a suitable $\wt{\C}$-linear topology and discuss some examples concerning Colombeau algebras of regular generalized functions.

First of all since $\G^\infty_E\subseteq \G_E$ the sharp topology induced by $\G_E$ on $\G^\infty_E$ gives the structure of a separated locally convex topological $\wt{\C}$-module to $\G^\infty_E$. The ultra-pseudo-seminorms obtained in this way are the original $\mP_i$ on $\G_E$ restricted to $\G^\infty_E$.

The moderateness properties of $\M_E^\infty$ allow us to define the valuation $\val^\infty_E:\M_E^\infty\to (-\infty,+\infty]$ as 
\[
\val^\infty_E((u_\eps)_\eps)= \sup\{b\in\R:\ \forall i\in I\ p_i(u_\eps)=O(\eps^b)\ \text{as}\ \eps\to 0\}
\]
which can be obviously extended to $\G^\infty_E$. This yields the existence of the ultra-pseudo-norm
\[
\mP_E^\infty(u):=\esp^{-\val^\infty_E(u)}
\]
on $\G_E^\infty$. Since for all $i\in I$ and $u\in\G^\infty_E$, $\val_{p_i}(u)\ge \val^\infty_E(u)$, the topology $\tau^\infty_\sharp$ determined by $\mP^\infty_E$ on $\G^\infty_E$ is finer than the topology induced by $\G_E$.

Adapting Proposition \ref{complete_G_E} to this situation one easily shows that when $E$ is a locally convex topological vector space with a countable base of neighborhoods of the origin, $\G^\infty_E$ with the topology $\tau^\infty_\sharp$ is complete. In fact assuming that $E$ is topologized through an increasing sequence of seminorms $\{p_k\}_{k\in\N}$, if $(u_n)_n$ is a Cauchy sequence in $\G^\infty_E$ then we can extract a subsequence $(u_{n_k})_k$ such that $\val^\infty_E(u_{n_{k+1}}-u_{n_k})>k$. This means that $\val_{p_i}(u_{n_{k+1}}-u_{n_k})>k$ for all $i\in \N$ and that there exists a decreasing sequence $\eps_k\searrow 0$, $\eps_k\le 2^{-k}$ such that $p_k(u_{n_{k+1},\eps}-u_{n_k,\eps})\le\eps^k$ for all $\eps\in(0,\eps_k)$. Defining $(h_{k,\eps})_\eps$ as in the proof of Proposition \ref{complete_G_E}, $(h_{k,\eps})_\eps\in\M^\infty_E$ and for all $k'\le k$ we have that $p_{k'}(h_{k,\eps})\le\eps^k$ on the interval $(0,1]$. As a consequence $u_\eps=u_{n_0,\eps}+\sum_{k=0}^\infty h_{k,\eps}$ is an element of $\M^\infty_E$ and for all $\overline{k}\ge 1$, $k\le\overline{k}$ and $\eps\in(0,\eps_{\overline{k}-1})$
\[
p_k(u_{n_{\overline{k}},\eps}-u_\eps)\le p_{\overline{k}}(u_{n_{\overline{k}},\eps}-u_\eps) \le c\eps^{\overline{k}-1}.
\]
When $k>\overline{k}$ we may write
\[
p_k(u_{n_{\overline{k}},\eps}-u_\eps)\le p_k(u_{n_{\overline{k}},\eps}-u_{n_k,\eps})+p_k(u_{n_{k,\eps}}-u_\eps)
\]
where as before $p_k(u_{n_{k,\eps}}-u_\eps)=O(\eps^{k-1})=O(\eps^{\overline{k}-1})$ and $p_k(u_{n_{\overline{k}},\eps}-u_{n_k,\eps})=O(\eps^{\overline{k}})=O(\eps^{\overline{k}-1})$ using the assumption $\val^\infty_E(u_{n_{k+1}}-u_{n_k})>k$ and a telescope sum argument. In this way we obtain that $\val^\infty_E(u_{n_{\overline{k}},\eps}-u_\eps)\ge \overline{k}-1$ and therefore $u_{n_k}\to u$ in $(\G^\infty_E,\tau^\infty_\sharp)$.\\
In conclusion we can say that if $E$ has a countable base of neighborhoods of the origin then the associated space $\G^\infty_E$ of regular generalized functions is a complete and ultra-pseudo-normed $\wt{\C}$-module.
\end{example}
\begin{example}
\label{example_GSinf}
\bf{The Colombeau algebra of $\S$-regular generalized functions}\rm 

A concrete example of $\G^\infty_E$ is given by the Colombeau algebra of $\S$-regular generalized functions $\GSinf(\R^n)$ introduced in \cite{Garetto:04, GGO:03}, whose definition is precisely $\G^\infty_E$ with $E=\S(\R^n)$. In this case we have that $\val^\infty_{\S(\R^n)}(u):=\sup\{b\in\R:\ \forall k\in\N\quad \sup_{x\in\R^n, |\alpha|\le k}(1+|x|)^k|\partial^\alpha u_\eps(x)|=O(\eps^b)\}$ and since $\mP^\infty_{\S(\R^n)}(uv)\le \mP^\infty_{\S(\R^n)}(u)\mP^\infty_{\S(\R^n)}(v)$ it turns out that $\GSinf(\R^n)$ is a topological $\wt{\C}$-algebra and a complete ultra-pseudo-normed $\wt{\C}$-module.
\end{example}
\begin{example}
\label{example_Ginf}
\bf{The Colombeau algebra $\Ginf(\Om)$}\rm

We recall that the \emph{Colombeau algebra of regular generalized functions} is the set $\Ginf(\Om)$ of all $u\in\G(\Om)$ having a representative $(u_\eps)_\eps$ satisfying the following condition
\beq
\label{estimate_regular}
\forall K\Subset\Om\ \exists N\in\N\ \forall\alpha\in\N^n\qquad \sup_{x\in K}|\partial^\alpha u_\eps(x)|=O(\eps^{-N})\ \text{as}\ \eps\to 0.
\eeq
\eqref{estimate_regular} defines the subset $\EMinf(\Om)$ of the set of moderate nets $\EM(\Om)$ and determines $\Ginf(\Om)$ as the factor $\EMinf(\Om)/\Neg(\Om)$. $\Ginf(\Om)$ can be seen as the intersection $\displaystyle\cap_{K\Subset\Omega}\Ginf(K)$ where $\Ginf(K)$ is the space of all $u\in\G(\Om)$ such that there exists a representative $(u_\eps)_\eps$ satisfying the condition
\beq
\label{estimate_GinfK}
\exists N\in\N\ \forall\alpha\in\N^n\qquad \sup_{x\in K}|\partial^\alpha u_\eps(x)|=O(\eps^{-N})\quad \text{as}\ \eps\to 0.
\eeq
Let us choose an exhausting sequence $K=K_0\subset K_1\subset K_2...$ of compact subsets of $\Omega$. We equip $\Ginf(K)$ with the locally convex $\wt{\C}$-linear topology determined by the usual ultra-pseudo-seminorms $\mP_i(u)=\esp^{-\val_{p_i}(u)}$ on $\G(\Om)$ where $p_i(u_\eps)=\sup_{x\in K_i, |\alpha|\le i}|\partial^\alpha u_\eps(x)|$, and the $\Ginf(K)$-ultra-pseudo-seminorm $\mP_{\Ginf(K)}$ defined via the valuation 
\beq
\label{val_Ginf(K)}
\val_{\Ginf(K)}(u)=\sup\{b\in\R:\ \forall\alpha\in\N^n\qquad \sup_{x\in K}|\partial^\alpha u_\eps(x)|=O(\eps^b)\}.
\eeq
It is clear that with this topology $\Ginf(K)$ is separated and metrizable. 

We want to prove that $\Ginf(K)$ is complete. As in the proof of Proposition \ref{complete_G_E} and the reasoning concerning $\Ginf_E$ in Example \ref{example_regular}, if $(u_n)_n$ is a Cauchy sequence in $\Ginf(K)$ we can extract a subsequence $(u_{n_j})_j$ and a sequence $\eps_j\searrow 0$, $\eps_j\le 2^{-j}$ such that for all $j\in\N$, $\val_{\Ginf(K)}(u_{n_{j+1},\eps}-u_{n_j,\eps})>j$ and $p_j(u_{n_{j+1},\eps}-u_{n_j,\eps})\le\eps^j$ on $(0,\eps_j)$. Define the net $(h_{j,\eps})_\eps$ as $u_{n_{j+1},\eps}-u_{n_j,\eps}$ on the interval $(0,\eps_j)$ and $0$ outside. By construction $(h_{j,\eps})_\eps$ satisfies \eqref{estimate_GinfK} and by Proposition \ref{complete_G_E} $u_\eps:=u_{n_0,\eps}+\sum_{j=0}^\infty h_{j,\eps}$ belongs to $\EM(\Om)$. More precisely for all $\alpha\in\N^n$
\[
\sup_{x\in K}|\partial^\alpha u_\eps(x)|\le \sup_{x\in K}|\partial^\alpha u_{n_0,\eps}(x)|+\sum_{j=0}^{|\alpha|}\sup_{x\in K, |\beta|\le |\alpha|}|\partial^\beta h_{j,\eps}(x)|+\sum_{j=|\alpha|+1}^\infty p_j(h_{j,\eps})
\]
\[
\le\sup_{x\in K}|\partial^\alpha u_{n_0,\eps}(x)|+\sum_{j=0}^{|\alpha|}\sup_{x\in K, |\beta|\le |\alpha|}|\partial^\beta h_{j,\eps}(x)|+\sum_{j=|\alpha|+1}^\infty\frac{1}{2^j},
\]
where $\sup_{x\in K, |\beta|\le |\alpha|}|\partial^\beta h_{j,\eps}(x)|=O(1)$ for all $j$ and $\alpha$. It follows that $(u_\eps)_\eps+\Neg(\Om)\in\Ginf(K)$ and adapting the estimates in the proof of Proposition \ref{complete_G_E} to our situation we obtain that for all $\overline{\jmath}\ge 1$
\beq
\label{est1}
p_{\overline{\jmath}}(u_{n_{\overline{\jmath}},\eps}-u_\eps)=O(\eps^{\overline{\jmath}-1})\quad \text{as}\ \eps\to 0.
\eeq
As a consequence $\sup_{x\in K, |\alpha|\le \overline{\jmath}}|\partial^\alpha(u_{n_{\overline{\jmath}},\eps}-u_\eps)(x)|=O(\eps^{\overline{\jmath}-1})$. If $j\ge\overline{\jmath}$ the assumption $\val_{\Ginf(K)}(u_{n_{j+1},\eps}-u_{n_j,\eps})>j$ leads to 
\beq
\label{est2}
\begin{split}
\sup_{x\in K, |\alpha|\le j}\ &|\partial^\alpha(u_{n_{\overline{\jmath}},\eps}-u_\eps)(x)|\\
&\le \sup_{x\in K, |\alpha|\le j}|\partial^\alpha(u_{n_{\overline{\jmath}},\eps}-u_{n_j,\eps)}(x)|+\sup_{x\in K_j,|\alpha|\le j}|\partial^\alpha(u_{n_j,\eps}-u_\eps)(x)|\\
&=O(\eps^{\overline{\jmath}})+O(\eps^{j-1})=O(\eps^{\overline{\jmath}-1}).
\end{split}
\eeq
\eqref{est2} shows that $\val_{\Ginf(K)}(u_{n_{\overline{\jmath}},\eps}-u_\eps)\ge\overline{\jmath}-1$ and combined with \eqref{est1} yields that $u_{n_j}$ is convergent to $u$ in $\Ginf(K)$. In this way $\Ginf(K)$ is a Fr\'echet $\wt{\C}$-module.

By definition $\Ginf(K')$ is a $\wt{\C}$-submodule of $\Ginf(K)$ when $K\subseteq K'$ and noting that $\mP_{\Ginf(K)}(u)\le\mP_{\Ginf(K')}(u)$ for all $u\in\Ginf(K')$, the topology on $\Ginf(K')$ is finer than the topology induced by $\Ginf(K)$ on $\Ginf(K')$. We are in the situation of Proposition \ref{prop_complete_initial}. Hence $\Ginf(\Om)$ equipped with the initial topology for the injections $\Ginf(\Om)\to\Ginf(K)$ is a complete locally convex topological $\wt{\C}$-module. More precisely since for every $\mP_i$ as above the estimate $\mP_i(u)\le\mP_{\Ginf(K_i)}(u)$ holds on $\Ginf(\Om)$, the initial topology on $\Ginf(\Om)$ is determined by a countable family of ultra-pseudo-seminorms and then $\Ginf(\Om)$ itself is a Fr\'echet $\wt{\C}$-module. This topology on $\Ginf(\Om)$ is finer than the sharp topology induced by $\G(\Om)$ on $\Ginf(\Om)$ and makes the multiplication of generalized functions in $\Ginf(\Om)$ continuous. 

Note that choosing $E=\E(\Om)$ in Example \ref{example_regular} we can construct the algebra $\G^\infty_{\E(\Om)}$. Obviously $\G^\infty_{\E(\Om)}\subseteq\Ginf(\Om)$ but they do not coincide since the estimates which concern the representatives in $\G^\infty_{\E(\Om)}$ require the same power of $\eps$ for all derivatives and all compact sets $K$. Finally $\mP_{\Ginf(K)}(u)\le\mP^\infty_{\E(\Om)}(u)$ for all $K\Subset\Om$ and $u\in\Ginf_{\E(\Om)}$.
\end{example}
\begin{example}
\label{example_Gcinf}
\bf{The Colombeau algebra $\Gcinf(\Om)$}\rm

$\Gcinf(\Om)$ denotes the algebra of generalized functions in $\Ginf(\Om)$ which have compact support. We want to endow this space with a $\wt{\C}$-linear topology. By the previous considerations each $\G^\infty_{\D_K(\Om)}$, $K\Subset\Om$, is a complete ultra-pseudo-normed $\wt{\C}$-module with valuation $$\val^\infty_{\D_K(\Om)}(u)=\sup\{b\in\R:\ \forall\alpha\in\N^n\ \sup_{x\in K}|\partial^\alpha u_\eps(x)|=O(\eps^b)\}.$$
Note that $\val_{\Ginf(K)}$ defined in \eqref{val_Ginf(K)} and $\val^\infty_{\D_K(\Om)}$ coincide on $\Ginf_{\D_K(\Om)}$. Repeating the reasoning of Example \ref{example_Gc}, $\G_K^\infty(\Om)$, the space of all generalized functions in $\Ginf(\Om)$ with support contained in $K$, is contained in $\Ginf_{\D_{K'}(\Om)}$, for every compact $K'$ containing $K$ in its interior. This inclusion is given by 
\[
\G^\infty_K(\Om)\to\Ginf_{\D_{K'}(\Om)}:u\to (u_\eps)_\eps+\Neg_{\D_{K'}(\Om)},
\]
where we choose a representative $(u_\eps)_\eps$ of $u$ in $\EMinf(\Om)\cap\M_{\D_{K'}(\Om)}$. It is clear that $\Ginf_{\D_{K'}(\Om)}$ is contained in $\Ginf(\Om)$. Since for every $K'_1,K'_2\Subset\Omega$ with $K\subseteq\text{int}(K'_1)\cap\text{int}(K'_2)$ and $u\in\Ginf_K(\Om)$ we have that $\val^\infty_{\D_{K'_2}(\Om)}(u)=\val^\infty_{\D_{K'_1}(\Om)}(u)$, we may define the valuation $\val^\infty_K$ on $\Ginf_K(\Om)$ as $\val^\infty_K(u)= \val^\infty_{\D_{K'}(\Om)}(u)$ where $K\subseteq\text{int}(K')\Subset\Omega$. In this way we can equip $\G^\infty_K(\Om)$ with the $\wt{\C}$-linear topology determined by the ultra-pseudo-norm $\mP_{\Ginf_K(\Om)}(u)=\esp^{-\val^\infty_K(u)}$. Note that this topology is finer than the one induced by $\G_K(\Om)$ on $\Ginf_K(\Om)$. 

Since the topology considered on $\Ginf_{\D_{K'}(\Om)}$ is finer than the topology induced by $\G_{\D_{K'}}(\Om)$ on $\Ginf_{\D_{K'}(\Om)}$ and $\G_K(\Om)$ is a closed subset of $\G_{\D_{K'}(\Om)}$ we have that $\Ginf_{K}(\Om)$ is closed in $\Ginf_{\D_{K'}(\Om)}$. Hence $\Ginf_K(\Om)$
is complete for the topology defined by $\mP_{\Ginf_K(\Om)}$. In analogy with the non-regular context examined before, $\Ginf_{K_2}(\Om)$ induces on $\Ginf_{K_1}(\Om)$ the original topology if $K_1\subseteq K_2$ and $\Ginf_{K_1}(\Om)$ is closed in $\Ginf_{K_2}(\Om)$.

At this point given an exhausting sequence $K_0\subset K_1\subset K_2....$ of compact sets, the strict inductive limit procedure equips $\Gcinf(\Om)=\cup_{n\in\N}\Ginf_{K_n}(\Om)$ with a complete and separated locally convex $\wt{\C}$-linear topology. Denoting the topologies on $\G_{K_n}(\Om)$ and $\Ginf_{K_n}(\Om)$ by $\tau_n$ and $\tau^\infty_n$ respectively and the inductive limit topologies on $\G_c(\Om)$ and $\Gcinf(\Om)$ by $\tau$ and $\tau^\infty$ respectively, we obtain that $\tau^\infty$ is the finest locally convex $\wt{\C}$-linear topology such that the embeddings $(\Ginf_{K_n}(\Om),\tau^\infty_n)\to(\Gcinf(\Om),\tau^\infty)$ are continuous. Moreover since the embedding maps $(\Ginf_{K_n}(\Om),\tau^\infty_n)\to(\G_{K_n}(\Om),\tau_n)\to(\Gc(\Om),\tau)$ are continuous, the topology $\tau^\infty$ on $\Gcinf(\Om)$ is finer than the topology induced by $\Gc(\Om)$.
\end{example}

\subsection{Continuity of a $\wt{\C}$-linear map $T:\G_E\to\G$}
We already argued on the continuity of a $\wt{\C}$-linear map between locally convex topological $\wt{\C}$-modules in Subsection 1.2, Theorem \ref{theorem_sem} and Corollary \ref{corollary_linear}. Here we focus our attention on $\wt{\C}$-linear maps where at least the domain is of a space of generalized functions $\G_E$ over $E$. In particular, we investigate the relationships between $\wt{\C}$-linearity and continuity with respect to the sharp topology by means of some examples. Before proceeding we recall that given locally convex topological vector spaces $(E,\{p_i\}_{i\in I})$ and $(F,\{q_j\}_{j\in J})$, by \eqref{est_gen_lin} a $\wt{\C}$-linear map $T:\G_E\to \G_F$ is continuous for the corresponding sharp topologies (``is sharp continuous''  for short) if and only if for all $j\in J$
\[
\mQ_j(Tu)\le C\max_{i\in I_0}\mP_i(u)
\]
for some finite subset $I_0$ of $I$, or in terms of valuations $\val_{q_j}(Tu)\ge -\log C + \min_{i\in I_0}\val_{p_i}(u)$.
In the following remark we discuss some examples of $\wt{\C}$-linear maps which are continuous.
\begin{remark}
\leavevmode
\label{rem_sharp_cont}
\begin{trivlist}
\item[(i)] When $E$ is a normed space with dim\,$E=n<\infty$, every $\wt{\C}$-linear map $T$ from $\G_E$ into a locally convex topological $\wt{\C}$-module $\G$ is continuous. 
\item[(ii)] We say that a $\wt{\C}$-linear map $T:\G_E\to\G_F$, where $E$ and $F$ are locally convex topological vector spaces, has a representative $t:E\to F$ if $t$ maps moderate nets into moderate nets and negligible nets into negligible nets, i.e. $(u_\eps)_\eps\in\M_E$ implies $(tu_\eps)_\eps\in\M_F$ and $(u_\eps)_\eps\in\Neg_E$ implies $(tu_\eps)_\eps\in\Neg_F$, and $Tu=[(tu_\eps)_\eps]$ for all $u\in\G_E$.

Any linear and continuous map $t:E\to F$ defines the $\wt{\C}$-linear and sharp continuous map $T:\G_E\to\G_F:u\to [(tu_\eps)_\eps]$. In fact, since $t$ is continuous $T$ is well-defined and sharp continuous and finally the $\C$-linearity of $t$ yields the $\wt{\C}$-linearity of $T$. 
\item[(iii)] The sharp continuity of a $\wt{\C}$-linear map $T$ with representative $t:E\to\C$ does not guarantee the continuity of the representative. Let $(E,\{p_i\}_{i\in I})$ be a non-bornological locally convex topological vector space and $t:E\to\C$ a linear bounded map which is not continuous. To provide an example of such a map $T$ we consider the space of regular generalized functions based on $E$ and we slightly modify the corresponding definition by requiring uniform estimates in the interval $(0,1]$. In other words we introduce $\underline{\G}^\infty_{\, E}\subseteq\Ginf_E$  by factorizing 
$\{(u_\eps)_\eps\in E^{(0,1]}:\ \exists N\in\N\, \forall i\in I\  \sup_{\eps\in(0,1]}\eps^{N}p_i(u_\eps)<\infty\}$ with respect to 
$\{(u_\eps)_\eps\in E^{(0,1]}:\ \forall i\in I\, \forall q\in\N\ \sup_{\eps\in(0,1]}\eps^{-q}p_i(u_\eps)<\infty\}$.
If $(u_\eps)_\eps$ is a representative of $u\in\underline{\G}^\infty_{\, E}$ then for some $N\in\N$ the set $\{\eps^{N}u_\eps,\ \eps\in(0,1]\}$ is bounded in $E$ and from the boundedness of $t$ we have that $|t(\eps^Nu_\eps)|=\eps^N|t(u_\eps)|\le c$. This means that $T:\underline{\G}^\infty_{\, E}\to\wt{\C}:u\to[(tu_\eps)_\eps]$ is well-defined and sharp continuous since $\val_{\wt{\C}}(Tu)\ge\val^\infty_E(u)$. \end{trivlist}
\end{remark}
Take now a pairing of vector spaces $(E,F,\bil)$ and endow $E$ with the weak topology $\sigma(E,F)$. It is clear that for each $y\in F$ the map $E\to \C:x\to \bil(x,y)$ is continuous and from $(ii)$ in Remark \ref{rem_sharp_cont} $\bil(\cdot,y):\G_E\to\wt{\C}:u\to \bil(u,y):=[(\bil(u_\eps,y))_\eps]$ is sharp continuous.
\begin{proposition}
\label{duality}
Let $(E,F,\bil)$ be a pairing and let $E$ be equipped with the weak topology $\sigma(E,F)$. If $T:\G_E\to\wt{\C}$ is a sharp continuous $\wt{\C}$-linear map with a representative $t:E\to\C$ then there exists $y\in F$ such that $Tu=\bil(u,y)$ for all $u\in\G_E$.
\end{proposition}
The proof of Proposition \ref{duality} needs a preparatory lemma.
\begin{lemma}
\label{lemma_kernel}
Under the assumptions of Proposition \ref{duality} on $E$ and $F$, if $T:\G_E\to\wt{\C}$ is a $\wt{\C}$-linear and sharp continuous map then there exists $\{y_i\}_{i=1}^N\subseteq F$ such that 
\beq
\label{inclusion}
\bigcap_{i=1}^N {\rm{ker}}\,\bil(\cdot,y_i) \subseteq {\rm{ker}}\,T.
\eeq
\end{lemma}
\begin{proof}
Since $T$ is continuous at the origin, there exists $\{y_i\}_{i=1}^N\subseteq F$ and $\eta>0$ such that $\max_{i=1}^N\vert\bil(u,y_i)\vert_\esp \le \eta$ implies $\vert Tu\vert_\esp \le 1$. Now if $u\in\bigcap_{i=1}^N {\rm{ker}}\,\bil(\cdot,y_i)$ then $\vert \bil(u,y_i)\vert_\esp=0$ for all $i=1,...,N$ and the same holds for $w=[(\eps^{-a})_\eps]u$ where $a$ is an arbitrary real number. Thus $\vert Tw\vert_\esp=\esp^{a}\vert Tu\vert_\esp\le 1$ which implies $\vert Tu\vert_\esp\le \esp^{-a}$ for all $a$ and therefore $u\in{\rm{ker}}\,T$.
\end{proof}
\begin{proof}[Proof of Proposition \ref{duality}]
By Lemma \ref{lemma_kernel} we know that there exists a finite number of $y_1,y_2,...,y_N$ in $F$ such that \eqref{inclusion} holds. Let $L$ be the $\wt{\C}$-linear map from $\G_E$ into ${\wt{\C}}^N$ given by $Lu=(\bil(u,{y_1}),\bil(u,{y_2}),...,\bil(u,{y_N}))$. The inclusion \eqref{inclusion} allows us to define $S:L(\G_E)\to\wt{\C}:(\bil(u,{y_1}),\bil(u,{y_2}),...,\bil(u,{y_N}))\to Tu$. Consider now the subset $V:=\{(\bil(u,{y_1}),\bil(u,{y_2}),...,\bil(u,{y_N})),\ u\in E\}$ of $\C^N$. By \eqref{inclusion} the map $s:V\to\C:(\bil(u,{y_1}),\bil(u,{y_2}),...,\bil(u,{y_N}))\to tu$ is a well-defined representative of $S$ and it can be obviously extended to a linear map $s':\C^N\to\C$. This means that if $(e_1,e_2,...,e_N)$ is the canonical basis of $\C^N$ and $s'(e_i)=\lambda_i$, $i=1,...,N$ then we have, for all $u\in\G_E$, that
\[
\begin{split}
Tu &= S(\bil(u,{y_1}),\bil(u,{y_2}),...,\bil(u,{y_N})) = [(s(\bil(u_\eps,{y_1}),\bil(u_\eps,{y_2}),...\bil(u_\eps,{y_N})))_\eps]\\ 
&= \biggl[\biggl(\sum_{i=1}^N\lambda_i\bil(u_\eps,{y_i})\biggr)_\eps\biggr] =\biggl[\biggl(\bil(u_\eps, \sum_{i=1}^N\lambda_i{y_i})\biggr)_\eps\biggr]= \bil\biggl(u, \sum_{i=1}^N \lambda_i y_i\biggr),
\end{split}
\]
where $\sum_{i=1}^N \lambda_i y_i \in F$.
\end{proof} 
In conclusion Proposition \ref{duality} combined with $(ii)$ in Remark \ref{rem_sharp_cont} and the classical results on pairings and continuity leads to the following statement: under the assumptions of Proposition \ref{duality}, $T:\G_E\to\wt{\C}$ is sharp continuous if and only if there exists $y\in F$ such that $T=\bil(\cdot,y)$ if and only if the representative $t:E\to\C$ is continuous.

\subsection{The topological dual $\L(\G_E,\wt{\C})$ when $E$ is a normed space}
The space of generalized functions $\G_E$ has a simple and interesting topological structure when $E$ is a normed space and we consider the ultra-pseudo-norm $\Vert u\Vert_{\G_E}:=\esp^{-\val_{\Vert\cdot\Vert_E}(u)}$. In Section 2 we equipped $\L(\G_E,\wt{\C})$ with three topologies: the weak topology $\sigma(\L(\G_E,\wt{\C}),\G_E)$, the strong topology $\beta(\L(\G_E,\wt{\C}),\G_E)$ and the polar topology $\beta_b(\L(\G_E,\wt{\C}),\G_E)$. The ultra-pseudo-norm introduced on $\G_E$ defines, as in the classical theory of normed spaces, a corresponding ultra-pseudo-norm on $\L(\G_E,\wt{\C})$ adding another $\wt{\C}$-linear topology to the list above. For the sake of generality we begin to discuss this topic in the context of an ultra-pseudo-normed $\wt{\C}$-module $\G$.
\begin{proposition}
\label{prop_dual_norm}
Let $\G$ be a topological $\wt{\C}$-module with topology determined by an ultra-pseudo-norm $\mP$. The map $\mP_{\L(\G,\wt{\C})}$ defined on $\L(\G,\wt{\C})$ by
\beq
\label{formula_norm}
\mP_{\L(\G,\wt{\C})}(T)=\inf\{C>0: \forall u\in\G\ \ \vert Tu\vert_\esp\le C\,\mP(u)\ \}
\eeq
is an ultra-pseudo-norm on $\L(\G,\wt{\C})$ and it coincides with $\displaystyle\sup_{\mP(u)=1}\vert Tu\vert_\esp$.
\end{proposition}
\begin{proof}
Since it is immediate, we do not prove that \eqref{formula_norm} has the properties which characterize an ultra-pseudo-norm. 
We note that when $u_0\neq 0$ in $\G$ the element $[(\eps^{\log\mP(u_0)})_\eps]u_0$ belongs to the set of $U:=\{u\in\G:\ \mP(u)=1\}$. Hence 
\begin{multline*}
\sup_{\mP(u)=1}\vert Tu\vert_\esp =\inf\{C>0:\, \forall u\in U\ \ \vert Tu\vert_\esp \le C\}\\= \inf\{C>0:\, \forall u\in\G\ \ \vert Tu\vert_\esp \le C \mP(u)\}.
\end{multline*}
\end{proof}
\begin{remark}
Denoting the topology on $\L(\G,\wt{\C})$ obtained via $\mP_{\L(\G,\wt{\C})}$ by $\tau$ we can write the chain of relationships
\beq
\label{chain1}
\sigma(\L(\G,\wt{\C}),\G)\,\preceq\, \tau\, \preceq\, \beta_b(\L(\G,\wt{\C}),\G)\,\preceq\, \beta(\L(\G,\wt{\C}),\G),
\eeq  
where $\preceq$ stands for ``is coarser than''.
\end{remark}
As in the theory of normed spaces we state the following result of completeness. The proof can be obtained by transferring the arguments of Proposition 3 in \cite[Chapter 3]{Royden:88} into the framework of ultra-pseudo-normed $\wt{\C}$-modules.
\begin{proposition}
\label{prop_complete_normed}
If $(\G,\mP)$ is an ultra-pseudo-normed $\wt{\C}$-module then the dual $(\L(\G,\wt{\C}),\mP_{\L(\G,\wt{\C})})$ is complete.
\end{proposition}
Before proceeding with $\G_E$ and the topological features of its dual $\L(\G_E,\wt{\C})$ we present an easy adaptation of the Banach-Steinhaus theorem to complete ultra-pseudo-normed $\wt{\C}$-modules. As observed in Subsection 1.2 every complete ultra-pseudo-normed $\wt{\C}$-module $\G$ is a complete metric space with metric $d(u_1,u_2)=\mP(u_1-u_2)$ and therefore a Baire space. This means that if $\G$ may be written as a countable union of closed subsets $S_n$ then at least one $S_n$ has nonempty interior. This fact allows us to prove the $\wt{\C}$-modules version of Osgood's theorem and the Banach-Steinhaus theorem.
\begin{theorem}
\label{Osgood}
Let $(\G,\mP)$ be a complete ultra-pseudo-normed $\wt{\C}$-module and $(T_i)_{i\in I}$ be a family of continuous functions defined on $\G$ with values in $\wt{\C}$. Suppose that for each $u\in\G$ the family $(T_i(u))_{i\in I}$ is bounded in $\wt{\C}$. Then there exist $u_0\in\G$, $\eta>0$ and $C>0$ such that $\vert T_i(u)\vert_\esp\le C$ for all $i\in I$ and $u\in\G$ with $\mP(u-u_0)\le\eta$.
\end{theorem}
\begin{theorem}
\label{Banach_Steinhaus}
Let $(\G,\mP)$ be a complete ultra-pseudo-normed $\wt{\C}$-module and $(T_i)_{i\in I}$ be a family of functions in $\L(\G,\wt{\C})$ such that $(T_i(u))_{i\in I}$ is bounded in $\wt{\C}$  for all $u\in\G$. Then there exists a constant $C>0$ such that $\Vert T_i\Vert_{\L(\G,\wt{\C})}\le C$ for all $i\in I$.
\end{theorem} 
\begin{proof}
By Theorem \ref{Osgood} we may find a set $B_{\eta}(u_0):=\{u\in\G:\ \mP(u-u_0)\le\eta\}$ and a constant $C>0$ such that $\vert T_i u\vert_\esp\le C$ for all $u\in B_\eta(u_0)$ and $i\in I$. Since $[(\eps^{\log(\mP(u)/\eta)})_\eps]u+u_0$ belongs to $B_\eta(u_0)$ when $u\neq 0$ in $\G$ it follows that
\begin{multline*}
\frac{\eta}{\mP(u)}\vert T_iu\vert_\esp=\vert T_i([(\eps^{\log(\mP(u)/\eta)})_\eps]u)\vert_\esp\\
\le\max\{\vert T_i([(\eps^{\log(\mP(u)/\eta)})_\eps]u+u_0)\vert_\esp,\vert T_iu_0\vert_\esp\}\le C.
\end{multline*}
Therefore $\vert T_iu\vert_\esp\le (C/\eta)\mP(u)$ for all $u\in\G$ and $i\in I$. This yields the uniform bound $\Vert T_i\Vert_{\L(\G,\wt{\C})}\le C/\eta $.
\end{proof}
For any normed space $E$ we already proved that $\L(\G_E,\wt{\C})$ is a complete ultra-pseudo-normed $\wt{\C}$-module for the ultra-pseudo-norm $\Vert\cdot\Vert_{\L(\G_E,\wt{\C})}$ defined by $\mP=\Vert\cdot\Vert_{\G_E}$ in \eqref{formula_norm}. The generalized functions belonging to $\G_{E'}$, when $E'$ is the topological dual of $E$ topologized through the norm $\Vert l\Vert_{E'}=\sup_{\Vert x\Vert_E=1}|l(x)|$, are particular elements of $\L(\G_E,\wt{\C})$.
\begin{proposition}
\label{prop_G_E'}
Let $E$ be a normed space. The map
\beq
\label{map_G_E'}
\G_{E'}\to \L(\G_E,\wt{\C}):v\ \to\ (u\to v(u):=[(v_\eps(u_\eps))_\eps])
\eeq
is a $\wt{\C}$-linear injection continuous with respect to $\Vert\cdot\Vert_{\G_{E'}}$ and $\Vert\cdot\Vert_{\L(\G_E,\wt{\C})}$.
\end{proposition}
\begin{proof}
First of all $|v_\eps(u_\eps)|\le \Vert v_\eps\Vert_{E'}\Vert u_\eps\Vert_E$ for all $(v_\eps)_\eps\in\M_{E'}$ and $(u_\eps)_\eps\in\M_E$. This implies that the map in \eqref{map_G_E'} is well-defined, $\wt{\C}$-linear and continuous by $\Vert v(\cdot)\Vert_{\L(\G_E,\wt{\C})}=\sup_{\Vert u\Vert_{\G_E}=1}\vert v(u)\vert_\esp\le \Vert v\Vert_{\G_{E'}}$. Concerning the injectivity, assume that $v(\cdot)=0$ in $\L(\G_E,\wt{\C})$ but $v\neq 0$ in $\G_{E'}$. This means that there exists a representative $(v_\eps)_\eps$ of $v$ such that $\Vert v_{\eps_n}\Vert_{E'}>\eps_n^q$ for some $q\in\N$ and a decreasing sequence $\eps_n\to 0$. Therefore, we may choose a sequence $(u_n)_n\subseteq E$ with $\Vert u_n\Vert_E=1$ for all $n$ such that $|v_{\eps_n}(u_{n})|>\eps_n^q$. Let now $(u_\eps)_\eps$ be the net in $E^{(0,1]}$ defined by $0$ on $(\eps_0,1]$ and $u_{n}$ on $(\eps_{n+1},\eps_n]$. Clearly $(u_\eps)_\eps\in\M_E$ and by construction $(v_\eps(u_\eps))_	\eps\not\in\Neg$. Therefore $v(u)\neq 0$ for $u=[(u_\eps)_\eps]\in\G_E$ which contradicts our hypothesis.
\end{proof}
In analogy with the Hahn-Banach theorem for normed spaces we may construct an element of the dual $\L(\G_E,\wt{\C})$ having an assigned value on some $u\in\G_E$. In the sequel we denote the complex generalized number $[(\Vert u_\eps\Vert_E)_\eps]$ by $\Vert u\Vert_E$.
\begin{proposition}
\label{prop_Hahn_Banach}
For any $u\in\G_E$ there exists $v\in\G_{E'}$ such that $\Vert v\Vert_{\G_{E'}}=1$ and $v(u)=\Vert u\Vert_E$.
\end{proposition}
\begin{proof}
Take a representative $(u_\eps)_\eps$ of $u$. By the Hahn-Banach theorem we have that for all $\eps\in(0,1]$ there exists $v_\eps\in E'$ such that $v_\eps(u_\eps)=\Vert u_\eps\Vert_E$ and $\Vert v_\eps\Vert_{E'}=1$. Hence $(v_\eps)_\eps\in\M_{E'}$ and  $v=[(v_\eps)_\eps]\in\G_{E'}$ satisfies the assertion.
\end{proof}
\begin{corollary}
\label{corollary_Hahn_Banach}
For all $u\in\G_E$
\beq
\label{norm_E}
\Vert u\Vert_{\G_E}=\sup_{\substack{v\in\G_{E'},\\ \Vert v\Vert_{\G_{E'}}=1}}\hskip-5pt\vert v(u)\vert_\esp .
\eeq
\end{corollary}
\begin{proof}
The right-hand side of \eqref{norm_E} is smaller than the left-hand side since the estimate $\vert v(u)\vert_\esp\le \Vert v\Vert_{\G_{E'}}\Vert u\Vert_{\G_E}$ holds for all $u\in\G_E$ and $v\in\G_{E'}$. By Proposition \ref{prop_Hahn_Banach} there exists $v\in\G_{E'}$ with ultra-pseudo-norm $1$ such that $v(u)=\Vert u\Vert_E$. Then $\vert v(u)\vert_\esp=\Vert u\Vert_{\G_E}$ and the equality in \eqref{norm_E} is attained.
\end{proof}
\begin{remark}
\label{remark_isometry}
Corollary \ref{corollary_Hahn_Banach} says that $(\G_E,\Vert\cdot\Vert_{\G_E})$ is isometrically contained in $(\L(\G_{E'},\wt{\C}),\Vert\cdot\Vert_{\L(\G_{E'},\wt{\C})})$. In particular applying this result to $\G_F$ with $F=E'$ we have that $(\G_{E'},\Vert\cdot\Vert_{\G_{E'}})$ is isometrically contained in $(\L(\G_{E},\wt{\C}),\Vert\cdot\Vert_{\L(\G_{E},\wt{\C})})$ when $E$ is reflexive.
\end{remark}
We conclude the paper proving the following proposition on bounded subsets.
\begin{proposition}
\label{prop_isometry}
Let $E$ be a normed space. $A\subseteq\G_E$ is $\Vert\cdot\Vert_{\G_E}$-bounded if and only if it is $\sigma(\G_E,\L(\G_E,\wt{\C}))$-bounded.
\end{proposition}
\begin{proof}
If $A$ is $\Vert\cdot\Vert_{\G_E}$-bounded then it is $\sigma(\G_E,\L(\G_E,\wt{\C}))$-bounded since the topology $\sigma(\G_E,\L(\G_E,\wt{\C}))$ is coarser than the sharp topology on $\G_E$ defined by $\Vert\cdot\Vert_{\G_E}$. Assume now that $A$ is $\sigma(\G_E,\L(\G_E,\wt{\C}))$-bounded. For all $v\in\G_{E'}\subseteq\L(\G_E,\wt{\C})$ we have that $\sup_{u\in A}\vert v(u)\vert_\esp <+\infty$. We can interpret the set $A$ as a family of maps in $\L(\G_{E'},\wt{\C})$ such that for all $v\in\G_{E'}$ the family $(v(u))_{u\in A}$ is bounded in $\wt{\C}$. Since $\G_{E'}$ is complete, by Theorem \ref{Banach_Steinhaus} there exists $C>0$ such that $\Vert u\Vert_{\L(\G_{E'},\wt{\C})}\le C$ for all $u\in A$. By Corollary \ref{corollary_Hahn_Banach} we conclude that $\sup_{u\in A}\Vert u\Vert_{\G_E}\le C$. This means that $A$ is $\Vert\cdot\Vert_{\G_E}$-bounded.
\end{proof}
As a consequence of Proposition \ref{prop_isometry} one may write 
\beq
\sigma(\L(\G_E,\wt{\C}),\G_E)\preceq\tau\preceq\beta_b(\L(\G_E,\wt{\C}),\G_E)=\beta(\L(\G_E,\wt{\C}),\G_E)
\eeq
where $\tau$ is the topology of the ultra-pseudo-norm $\Vert\cdot\Vert_{\L(\G_E,\wt{\C})}$ on $\L(\G_E,\wt{\C})$.

\bibliographystyle{abbrv}
\newcommand{\SortNoop}[1]{}

\end{document}